\documentclass[12pt]{amsart}
\usepackage[nohug]{diagrams}
\usepackage[mathscr]{eucal}
\usepackage{amscd}
\usepackage{amsfonts}
\usepackage{amsmath}
\usepackage{amsthm}
\usepackage{amssymb}
\usepackage{latexsym}
\usepackage[left=2cm,top=2cm,right=3cm,foot=2cm,headsep=1cm]{geometry}

\begin{document}

\newcommand{\nc}{\newcommand}
\newcommand{\bom}{{_{\mathbf{\omega}}}}
\newcommand{\st}{\divideontimes}
\def\neweq{\setcounter{theorem}{0}}
\newtheorem{theorem}[]{Theorem}
\newtheorem{proposition}[]{Proposition}
\newtheorem{corollary}[]{Corollary}
\newtheorem{lemma}[]{Lemma}
\theoremstyle{definition}
\newtheorem{definition}[]{Definition}
\newtheorem{remark}[]{Remark}
\newtheorem{conjecture}[equation]{Conjecture}
\newcommand{\dis}{{\displaystyle}}
\def\question{\noindent\textbf{Question.} }
\def\remark{\noindent\textbf{Remark.} }
\def\proof{\medskip\noindent {\textsl{Proof.} \ }}
\def\endproof{\hfill$\square$\medskip}
\def\str{\rule[-.2cm]{0cm}{0.7cm}}
\newcommand{\beq}{\begin{equation}\label}
\newcommand{\aand}{\quad{\text{\textsl{and}}\quad}}
%
%
\def\x{\stackrel{\rightarrow}{x}}
\newcommand{\iso}{{\;\;\stackrel{\sim}{\longrightarrow}\;\;}}
\newcommand{\cd}{\!\cdot\!}
\newcommand{\vi}{${\sf {(i)}}\;$}
\newcommand{\vii}{${\sf {(ii)}}\;$}
\newcommand{\viii}{${\sf {(iii)}}\;$}
\newcommand{\viv}{${\sf {(iv)}}\;$}
\newcommand{\sset}{\subset}
\newcommand{\G}{\Gamma}
\newcommand{\tdr}{{\widetilde{\sf{DR}}}}
\newcommand{\id}{{{\mathtt {Id}}}}
\newcommand{\im}{{{\mathtt {Im}}}}
\newcommand{\Ind}{{{\mathtt {Ind}}}}
\newcommand{\sq}{{$\enspace\square$}}
\newcommand{\into}{\,\,\hookrightarrow\,\,}
\newcommand{\too}{\,\,\longrightarrow\,\,}
\newcommand{\onto}{\,\,\twoheadrightarrow\,\,}
\newcommand{\bull}{^{^{_\bullet}}}
\newcommand{\sto}{\rightsquigarrow}
\newcommand{\ad}{{\mathtt{{ad}}}}
\newcommand{\Ad}{{\mathtt{{Ad}}}}
\newcommand{\Lie}{{\mathtt{Lie}^{\,}}}
\newcommand{\Spec}{{\mathtt{Spec}}}
\newcommand{\Hilb}{{\mathtt{Hilb}}}
\newcommand{\End}{{\mathtt{End}}}
\newcommand{\Endo}{{\mathtt{End}^{^{\mathsf{opp}}}}}
\newcommand{\Hom}{{\mathtt{Hom}}}
\newcommand{\Aut}{{\mathtt{Aut}}}
\newcommand{\DGAut}{{\mathtt{DGAut}}}
\newcommand{\Mod}{{\mathtt{Mod}}}
\newcommand{\MOD}{{\mathtt{MOD}}}
\newcommand{\Ideals}{{\mathtt{Ideals}}}
\newcommand{\DGMod}{{\mathtt{DGMod}}}
\newcommand{\Modi}{{\mathtt{Mod}}_{\infty}}
\newcommand{\Hi}{{\mathcal{H}}_{\infty}}
\newcommand{\Di}{{\mathcal{D}}_{\infty}}
\newcommand{\Com}{{\mathtt{Com}}}
\newcommand{\Ker}{{\mathtt{Ker}}}
\newcommand{\Coker}{{\mathtt{Coker}}}
\newcommand{\Hominf}{\raisebox{0ex}[2ex][0ex]
{$\overset{\,\infty}{\raisebox{0ex}[1ex][0ex]{$\mathtt{Hom}$}}$}}
\newcommand{\Span}{{\mathtt{span}}}
\newcommand{\Ann}{{\mathtt{Ann}}}
\renewcommand{\dim}{{\mathtt{dim}}}
\renewcommand{\deg}{{\mathtt{deg}}}
\newcommand{\Par}{{\mathtt{Part}}}
\newcommand{\co}{{\mathtt{cont}}}
\newcommand{\rk}{{\mathtt{rk}}}
\newcommand{\ch}{\chi}
\newcommand{\vr}{\boldsymbol{\varrho}}
\newcommand{\lrho}{\acute{\rho}}
\newcommand{\rrho}{\grave{\rho}}
\newcommand{\dr}{{{\mathsf {DR}}}}
\newcommand{\chom}{{\mathcal{H}om}}
\newcommand{\grd}{{\mathtt{gr}}}
\newcommand{\Ext}{{\mathtt{Ext}}}
\newcommand{\trace}{{\mathtt{Trace}}}
\newcommand{\Tr}{{{\mathtt{Tr}}}}
\newcommand{\tr}{{{\mathtt{Tr}}}}
\newcommand{\irr}{{{\mathsf {Irred}}}}
\newcommand{\rad}{^{^{\mathsf{rad}}}}
\newcommand{\hr}{{\mathfrak{h}^{^{\mathsf{reg}}}}}
\newcommand{\gr}{{\mathfrak{g}^{^{_{\mathsf{rs}}}}}}
\newcommand{\dd}{{\mathcal{D}}}
\newcommand{\GL}{\mathtt{GL}}
\newcommand{\SL}{\operatorname{SL}}
\newcommand{\SSp}{\mathtt{Sp}}
\newcommand{\hh}{{\mathsf{H}}}
\newcommand{\thh}{{\mathsf{\tilde{H}}}}
\newcommand{\ehe}{{\mathbf{e}\!\mathsf{H}\!\mathbf{e}}}
\newcommand{\eve}{{\mathbf{e}\mathcal{H}\mathbf{e}}}
\newcommand{\e}{{\mathbf{e}}}
\newcommand{\cc}{{\mathbf{c}}}
\newcommand{\ka}{\kappa}
\newcommand{\MM}{{\mathcal{M}}}
\newcommand{\RR}{{\mathcal{R}}}
\newcommand{\AB}{{\boldsymbol{A}}}
\newcommand{\BA}{{\boldsymbol{B}}}
\newcommand{\LB}{{\boldsymbol{L}}}
\newcommand{\sym}{{\mathsf{Sym}}}
\newcommand{\dertr}{{\mathtt{Der}}_{\!_{\mathsf {tr}}}}
\newcommand{\difp}{\mathtt{Diff}_{_{\!{\mathscr{P}}}}}
\newcommand{\derp}{\mathtt{Der}_{_{\!{\mathscr{P}}}}}
\newcommand{\sgn}{{\mathsf{sign}}}
\newcommand{\Ug}{{{\cal{U}}\g}}
\newcommand{\Uh}{{{\cal{U}}\h}}
\newcommand{\triv}{{\mathtt {triv}}}
\newcommand{\sign}{{\mathtt {sign}}}
\newcommand{\sch}{\boldsymbol{\mathsf{s}}}
\newcommand{\mon}{\boldsymbol{\mathsf{m}}}
\newcommand{\aut}{{\mathtt {Aut}}}
\newcommand{\tdash}{\mbox{\tiny{-}}}
\newcommand{\Th}{\Theta}
\newcommand{\lL}{\underround}
\newcommand{\rR}{{}}
\newcommand{\bphi}{\bar{\phi}}
\newcommand{\ee}{{\mathfrak{e}}}
\newcommand{\h}{{\mathfrak{h}}}
\newcommand{\cm}{{\mathbb{M}}}
\newcommand{\var}{\kappa}
\newcommand{\dvar}{{\mathsf{D}}_\kappa}
\newcommand{\Eu}{{{\mathtt {Eu}}}}
\newcommand{\om}{\Omega}
\newcommand{\A}{{\mathsf{A}}_{\infty}}
\newcommand{\Om}{\Omega}
\newcommand{\VV}{\mathsf{Vect}}
\newcommand{\wh}{\widehat}
\newcommand{\bh}{\bar{H}}
\newcommand{\la}{\label}
\newcommand{\greg}{{\mathfrak{g}^{^{_{\mathsf{reg}}}}}}
\newcommand{\hreg}{{\mathfrak{h}^{^{_{\mathsf{reg}}}}}}
\newcommand{\V}{{\mathscr{V}}_\var}
\newcommand{\PP}{{\mathscr{P}}}
\newcommand{\eij}{{\epsilon_{_{ij}}}}
\newcommand{\eps}{{\epsilon}}
\newcommand{\bul}{^\bullet\!}
\newcommand{\bulp}{^\bullet_{_{\!\PP\!}}}
\newcommand{\hv}{{\mathcal{H}}}
\newcommand{\va}{\varkappa}
\def\bs#1{\boldsymbol{#1}}
\def\ms#1{\mathcal{#1}}
\def\C{{\mathbb{C}}}
\def\X{{\bar{X}}}
\def\Y{{\bar{Y}}}
\def\j{{\bar{j}}}
\def\i{{\bar{i}}}
\def\g{{\bar{g}}}
\def\Q{{\mathbb{Q}}}
\def\R{{\mathsf{R}}}
\def\qq{{\mathcal{Q}}}
\def\ms#1{\mathcal{#1}}
\def\AA{{\mathbb{A}}}
\def\II{{\mathbb{I}}}
\def\rr{{\mathcal{R}}}
\def\Z{{\mathbb{Z}}}
\def\oo{{\mathcal O}}
\def\K{\boldsymbol{K}}
\def\KK{ \tilde{\boldsymbol{K}}}
\def\LL{ \tilde{\boldsymbol{L}}}
\def\L{\boldsymbol{L}}
\def\M{\boldsymbol{M}}
\def\N{\boldsymbol{N}}
\def\xx{\boldsymbol{x}}
\def\E{\boldsymbol{E}}
\def\U{{\mathcal U}}
\def\ccr{\C[R]^W_{_{\sf reg}}}
\def\ccz{\Z[R]^W}
\def\up{{\mathcal U}^{\PP\!}\!}
\def\upa{{{\mathcal U}^{\PP\!}\!A}}
\def\pp{_{_{\!\PP}}}
\def\ll{{\mathcal L}}
\def\LL{{\mathbf{L}}}
\def\D{{\mathcal{D}}}
\def\rep{{\mathsf{Rep}}}
\def\Pf{{\it Proof}}
\def\CP{{\mathbb{C}\mathbb{P}}}
\def\O{{\sf O}}
\def\aa{{\mathcal{A}}}
\def\eu{{\mathsf{eu}}}
\def\reg{{\!}^{^{\mathsf{reg}}}}
\def\bb{{\mathcal{B}}}
\def\T{{\mathsf T}}
\def\sll2{{\mathfrak{s}\mathfrak{l}}_2}
\def\gln{{\mathfrak{g}\mathfrak{l}}_n(\C)}
\def\gl{{\mathfrak{g}\mathfrak{l}}}
\def\sln{{\mathfrak{s}\mathfrak{l}}_n(\C)}
\def\glinfty{{\mathfrak{g}\mathfrak{l}}_\infty}
\def\Cr{\C^{\mathsf{reg}}}
\def\CC{{\mathcal{C}}}
\def\HH{{\mathcal{H}}}
\def\F{{\mathcal{F}}}
\def\FF{{\mathsf{F}}}
\def\s{{\mathbb{S}}}
\def\sn{{{\mathbb{S}}_n}}
\def\w{{\mathsf{W}}}
\def\bpa{{{\mathsf{B}}_{\bullet}^\PP\!A}}
\def\ccirc{{{}_{^{^\circ}}}}
\def\dertr{{\der_{\!_{\om\otimes\tr}}}}
\newcommand{\hdv}{{{\vphantom{H}
\smash{\buildrel_{\,_{_{_{\bullet\bullet}}}}\over {\sf H}}}}_{\!\kappa}}

\numberwithin{equation}{section}
\title{DG-models of projective modules and Nakajima quiver varieties}
\author{Farkhod Eshmatov}
\address{Department of Mathematics,
Cornell University, Ithaca, NY 14853-4201, USA}
\email{eshmatov@math.cornell.edu}
\maketitle
\section{ Introduction}
\la{Sect1}
This paper is inspired by the recent work of Berest and Chalykh \cite{BC}
on the right ideals of the first Weyl algebra $A_1(\mathbb C)$ and
Calogero-Moser spaces. The main result of \cite{BC} is an
explicit construction of the Calogero-Moser correspondence refining
the earlier work of Berest-Wilson \cite{BW1}.

The purpose of this paper is to extend the ideas and techniques of
\cite{BC} to a broader class of algebras of geometric origin. More
specifically, we will study right ideals in \textit{quantized}
coordinate rings of Kleinian singularities $\mathbb C^2 \!/ \Gamma$,
where $\Gamma$ is a finite cyclic subgroup of $SL_2 (\mathbb C)$.
For a fixed $\Gamma$, such rings form a family of noncommutative
algebras $O^ \tau$ (parametrized by the elements $\tau$ of the
group algebra $\mathbb C \Gamma$), whose properties are similar to
the properties of the Weyl algebra. Specifically, like $A_1 (\mathbb
C)$, in case of generic $\tau$'s the rings $O^\tau$ are simple,
hereditary, Noetherian domains, having no nontrivial
finite-dimensional representations. However, unlike $A_1$, they have
a nontrivial $K$-group. A conjectural description of stably free
ideals of $O^ \tau$\!, generalizing the work of Berest-Wilson, was
suggested by Crawley-Boevey and Holland (see \cite{CB}). 
Recently, Baranovsky, Ginzburg and Kuznetsov
\cite{BGK} refined and proved this conjecture using the methods of
noncommutative projective geometry. The main idea behind BGK's work,
which in the case of  $A_1$ was exploited earlier in \cite{LeB} and
\cite{BW2}, consists in replacing $O^ \tau$ by a graded algebra
$\mathbb B^\tau$\!, which, by analogy with geometric case, can be
treated as the homogeneous coordinate ring of a noncommutative projective
variety. Projective modules over $O^ \tau$ can then be extended to
certain ``vector bundles''    on such a ``variety'' and the latter can be
classified using the standard tools from algebraic geometry (the
Beilinson spectral sequence and Barth's monads). Despite its
naturality, this geometric approach has some disadvantages. First,
it is fairly complicated and far from being explicit. Second, it
involves a lot of choices (most notably the choice of filtration on
the given algebra $O^ \tau$), which are not intrinsic to the original
problem. Third, it hides some interesting ``affine''  features of the
problem, present in the case of the Weyl algebra: namely, the action
of the Dixmier automorphism group on the ideal classes and the
equivariance of the corresponding classifying map.

In the present paper we will give a new proof of the
Crawley-Boevey-Holland conjecture, which is free from the above
disadvantages. As in \cite{BC}, our construction is elementary and
independent of the choice of a filtration on $O^ \tau$; it leads to a
completely explicit description of ideals of $O ^ \tau$\!, and more
importantly, it is $G$-equivariant with respect to a certain ``large"
automorphism group $G$, which acts naturally on both the space of
ideal classes and the associated quiver varieties $\mathfrak M^{\tau}$.
This brings the picture with Kleinian singularities closer to the
original example of the Weyl algebra and raises many interesting
questions regarding the action of the group $G$ on $\mathfrak M^{\tau}$
(cf. \cite{BW1}).

\vspace{0.2cm}
\noindent
$\mathbf {Acknowledgments.}$ I would like to thank Oleg Chalykh and George Wilson for many
 useful discussions and valuable comments. I am especially grateful to my
advisor, Yuri Berest,  for his guidance and support  throughout this paper.

\section{Background and Statement of Results}
\subsection{ The algebras $B^\tau$  and $O^\tau$}
Let $(L, \omega)$ be a two-dimensional symplectic vector space
with symplectic form $\omega$, and let $\Gamma $ be a
finite subgroup of $Sp(L,\omega)$. We can extend the natural
(contragradient) action of $\Gamma$ on $L^*$ diagonally to
$TL^*$, the  tensor algebra of $L^*$ over $\mathbb C$, and
define $R$ to be the crossed product of $TL^*$ with 
$\Gamma$.
The form $\omega$ is a skew symmetric element of $L^* \otimes L^* \subset
TL^{*} \subset R $,
so for each $\tau \in Z(\mathbb C \Gamma)$  we can define 
\begin{eqnarray}
&& B^\tau = R  / R(\omega - \tau)R  , \nonumber \\
&& O^\tau = eB^\tau e , \nonumber
\end{eqnarray}
where $e$ is the symmetrizing idempotent $\sum_{g\in \Gamma} g / | \Gamma | 
$ in $\mathbb C\Gamma \subset B^\tau$ .
The algebras $B^\tau$ and $O^\tau$ have been introduced and studied by
W.Crawley-Boevey and M.Holland in \cite{CBH}.

It is convenient to choose a symplectic basis $\{e_ x , e_ y \}$ in $L$
and identify $L$ with $\mathbb C^2$, and $Sp(L,\omega)$ with $SL_2(\mathbb C)$.
If $\{ x,  y \}$ is the dual basis in $L^*$ then we have an algebra 
isomorphism:
\begin{equation}
\la{rel1}
 B^\tau \cong R / R( x y -  y  x -\tau)R \, , 
\end{equation}
where $ R \cong \mathbb C \langle x , y \rangle \ast \Gamma$ is a crossed
product of the free algebra on two generators with the group $\Gamma$.

In this paper we will be concerned with the case when $\Gamma $ is
a cyclic group $\mathbb Z_m$. One can give a more elementary description of 
$B^\tau$ in this case. We fix an embedding $\Gamma \hookrightarrow SL_2(\mathbb C)$
so that $L$ decomposes as  $\epsilon \oplus \epsilon ^{-1}$,
where $\epsilon$ is a primitive character of $\Gamma$.
Now we choose a basis $\{ x, y\}$ in $L^*$ so that $\Gamma$
acts on $x$ by $\epsilon$ and on $y$ by $\epsilon^{-1}$. Then
there is an algebra isomorphism from $B^\tau$ to the quotient
of $R$ by the following relations:
\begin{eqnarray}
\la{rel}
&& g \cdot x = \epsilon (g)\, x \cdot g \quad  , \quad  g \cdot y = \epsilon ^{-1} 
(g)\, y\cdot g ,
\quad  \forall g\in \Gamma, \\
&& x \cdot y - y \cdot x =\tau.
\end{eqnarray}
The corresponding algebras $O^\tau$ are called , in this case,
 \textit{the type A deformations of Kleinian singularities}.
They were studied earlier by Hodges \cite{H} and Smith \cite{Sm}.

The homological and ring-theoretical properties of $O^\tau$
depend drastically on values of the parameter $\tau$.
Through the McKay correspondence we can associate to the group $\Gamma$
the affine Dynkin graph of type $A$.
The group algebra $\mathbb C \Gamma$  is then identified
with the dual of the space spanned by the simple roots of the corresponding
affine root system and, following \cite{CBH}, we say that $\tau \in \mathbb C \Gamma$ 
is \textit{generic} if it does not belong to any root hyperplane in $\mathbb C \Gamma$.

From now on we will assume  $\tau$ to be generic. In this case
$B^\tau$ and $O^\tau$ are Morita equivalent, the equivalence $F :
Mod(B^\tau) \to  Mod(O^\tau) $ between the categories of right
modules is given by $ M \mapsto M \otimes_{B^\tau} B^\tau e$ (see
\cite{CBH}, Theorem $0.4$).
\subsection{Nakajima varieties }\label{SNv}
Given a pair  $ (U,W) $ of finite dimensional $ \Gamma $-modules 
consider
the space of $\Gamma$-equivariant linear maps
\begin{equation}
\la{Qv1} {\mathbb M}_{\Gamma}(U,W) = Hom_{\Gamma}(U, U \otimes L ) 
\oplus Hom_{\Gamma}(W,U)\oplus Hom_{\Gamma}(U,W) .
\end{equation}
The group $ G_{\Gamma}(U) $ of $ \Gamma $-equivariant automorphisms
of $ U $ acts on $ {\mathbb M}_{\Gamma}(U,W) $ in the natural way:
$$ g(B, \i, \j)=(g B g^{-1}, g \i, \j g^{-1}), $$
and this action is free on the subvariety
 $\tilde {\mathfrak M}^{\tau}_{\Gamma}(U,W) \subseteq \mathbb M_{\Gamma}(U,W) $ 
defined by the conditions 
\begin{eqnarray}
\la{NV}
&& \, (i)\, [B,B] + \tau \vert_{U} = \i \j  \, ,\\*[1ex] \nonumber
&& (ii)\,  \mbox{ There is no proper submodule } U' \subset U  \\*[1ex]
&& \mbox{such that } B(U')\subset U' \otimes L \mbox{ and }
\i (W) \subset U'.
  \nonumber
\end{eqnarray}
Here $[B, B]$ stands for the composition of the following maps
$$ U \overset B \longrightarrow U \otimes L
\stackrel{B \otimes id_{L}}{\longrightarrow}
U \otimes L \otimes L  \stackrel{id_{U} \otimes
\omega}{\longrightarrow} U\otimes \mathbb C \cong U \,  .$$
\begin{definition}
\la{D01} The Nakajima variety associated to the pair $ (U,W) $
is defined by 
$$ \mathfrak M^{\tau}_{\Gamma}(U,W):=  \tilde
{\mathfrak M}^{\tau}_{\Gamma}(U,W) /\!\!/ G_{\Gamma}(U) \, .$$
\end{definition}
The relation of this variety to the original definition of Nakajima \cite{N} 
can be obtained via the McKay correspondence (see e.g. [BGK]).

As $L\cong \epsilon \oplus\epsilon^{-1}$, we can write down the map
$B: U \to U \otimes L$ in the following form
\begin{equation}
\la{nk}
B(v) = \X(v) \otimes \epsilon + \Y(v) \otimes \epsilon^{-1}
\end{equation}
with $\X,\Y \in End_{\mathbb C}(U)$. The points of the Nakajima variety 
$ \mathfrak M^{\tau}_{\Gamma}(U,W)$
can be represented then by quadruples  $ (\X,\Y,\i,\j) \in End(U)^{\oplus 2} \oplus
 Hom_{\Gamma}(W,U) \oplus Hom_{\Gamma}(U,W) $ satisfying
\begin{equation}
\la{ma}
\X\,\Y - \Y\,\X + \bar T = \i\, \j\
\end{equation}
\begin{equation}
\la{mg}
\X\,\bar{G} =\epsilon(g)\, \bar{G}\,\X  \quad  , \quad  \Y\,\bar{G}
=\epsilon^{-1}(g)\, \bar{G}\,\Y
\end{equation}
where $\bar T$ and $ \bar{G} $ are the endomorphisms  corresponding the action of  $\tau$ and 
$ g\in \Gamma $ in $U$ respectively.

In case when $W$ is a one-dimensional $\Gamma$-module
with character $\chi_{W}$ we can think of $\i$ and $\j$ just as linear maps 
$\i \in Hom(W,U)$ and $\j \in Hom(U,W)$ satisfying the conditions:
\begin{equation}
\la{alt}
\i (w.g) = \chi_W(g)\, \i(w) \, , \, \j (v.g) = \chi_W(g)\, \j(v) \mbox{ 
for all }
g\in \Gamma \mbox{ and } w \in W , v\in U . 
\end{equation}
\subsection{ Statement of Results }
We start this section by reminding the reader the following result 
due to Berest and Wilson.
\begin{theorem}(\cite{BW1}, \cite{BW2})
\la{bw}
$(a)$ There is a natural bijection between the space $\mathcal R$ of isomorphism 
classes of right ideals of the Weyl algebra $A_1(\mathbb C)$ and the union 
$C=\bigsqcup C_n$ of Calogero-Moser algebraic varieties:
\begin{equation}
 C_n = \{ X, Y \in M_n(\mathbb C)\, \vert \, rk (XY - YX - Id) =1 \} / GL_n(\mathbb C) , \nonumber
\end{equation}
where $GL_n(\mathbb C)$ acts on $(X, Y)$ by simultaneous conjugation .\\
$(b)$ The automorphism group $G=Aut_{\mathbb C}(A_1)$ acts naturally on the varieties
$C_n$, and this action is transitive for each $n=0, 1, 2,...$\\
$(c)$ The bijection $\mathcal R \longleftrightarrow C$ is equivariant under $G$, and thus 
the varieties $C_n$ can be identified with the orbits of natural action of $G$ on set 
$\mathcal R$ of ideal classes.
 \end{theorem}
The algebras $O^\tau$ are obvious generalizations of $A_1$ and one might
expect that a result similar to Theorem \ref{bw} holds for the ideals of $O^\tau$.
 In fact, Crawley-Boevey and Holland have conjectured that there is a bijection 
between the space of isomorphism classes of ideals of $O^\tau$ and certain Nakajima 
varieties related to $\Gamma$. Such a classification of ideals in terms of ``Nakajima data'' 
suggests the existence of some finite dimensional module attached to each ideal.
The natural candidates for such modules would be finite dimensional representations
of the algebra $O^\tau$\!, but since $O^\tau$ is simple such representations do not exist. 
Nevertheless, stretching the notion of a module one may
overcome this problem. To be precise, we would like to extend the category
of modules over $O^\tau$ to the category of DG modules over a certain DG algebra
closely related to $O^\tau$. In this extended category we will construct objects 
whose isomorphism classes are in a natural bijection with
isomorphism classes of ideals in $O^\tau$, and from which we can extract
the Nakajima data corresponding to a given ideal. The idea of this approach 
goes back to \cite{BC}, where the ideals of $A_1(\mathbb C)$ are ``modelled''
by certain $\A$-modules. 

First, we would like to give a coarse, $K$-theoretic, classification 
of ideals of $O^\tau$.  Let $ K_0(\Gamma)$, $ K_0(B^\tau)$ and  
$ K_0(O^\tau)$ be the Grothendieck groups of the algebras $\mathbb C\Gamma$, $B^\tau$ and $O^\tau$ 
respectively. Then, by a well-known theorem of Quillen the induction functor $ P \mapsto 
P\otimes_{\mathbb C\Gamma} B^\tau$ gives an isomorphism of groups $K_0(\Gamma) \cong K_0(B^\tau)$.
Further, since $B^\tau$ and $O^\tau$ are Morita equivalent algebras,
the corresponding equivalence induces another isomorphism 
$K_0(B^\tau) \cong K_0(O^\tau)$. We will use these isomorphisms to identify 
$K_0(B^\tau)$ and $K_0(O^\tau)$ with $K_0(\Gamma)$. 
The map assigning to a finite dimensional module of $\Gamma$ its dimension 
extends to a group homomorphism $ \dim : K_0(\Gamma) \to \mathbb Z$.
 By Proposition \!\ref{neb} below, under our identification $K_0(O^\tau) \cong K_0(\Gamma)$ 
this homomorphism corresponds to the rank function on projective modules 
$ \mathtt{rk} : K_0(O^\tau) \to \mathbb Z$, and therefore each ideal of $O^\tau$ has 
dimension $1$ in $K_0(\Gamma)$. Now, according to \cite{BGK} Proposition 
\!$1.3.11$(see also Lemma \!\ref{P1} below) for each $[P]\in K_0(\Gamma)$ with $\dim \, [P] =1 $,
 there is a unique pair $(V,W) $ of finite-dimensional $\Gamma$-modules such 
that $\dim\, W =1$ , $V \nsubseteq \mathbb C\Gamma$ and 
$[P] = [W] + [V] \cdot ([L]-2[W_0])$, where $[W_0]\in K_0(\Gamma)$
is the class of the trivial representation of $\Gamma$. 

Given such $[P]\in K_0(\Gamma)$, we let
\begin{eqnarray}
&&\mathcal R (V,W):=\{ \mbox{Isomorphism classes of finitely generated projective $B^\tau$-modules }  M \\ 
&& \ \ \ \ \ \ \ \ \ \ \ \ \ \ \ \    \mbox{ such that } [M] = [P]  \mbox{ in } K_0(B^\tau)=
K_0(\Gamma) \} \nonumber\\
&&\mathcal R^{\prime} (V,W):=\{ \mbox{Isomorphism classes of ideals } N \mbox{ of }  O^\tau
\mbox{such that } [N] = [P] \mbox{ in } \\
&& \ \ \ \ \ \ \ \ \ \ \ \ \ \ \ \ K_0(O^\tau)=K_0(\Gamma) \}  \nonumber.
\end{eqnarray} 
Then the Morita equivalence $ F:\Mod (B^\tau) \to \Mod (O^\tau)$ 
induces a natural bijection
$$\mathcal R(V,W) \iso \mathcal R^\prime (V,W) \, , \, (M) \mapsto (Me)$$
Thus, the problem of classifying the ideals of $O^\tau$ is equivalent to 
classifying projective $B^\tau$-modules in
\begin{equation}
\la{eq1}
\mathcal R = \bigsqcup_{V,W} \mathcal R (V,W) \,  , 
\end{equation}
where $V$ runs over the set of isomorphism classes of all finite dimensional 
$\Gamma$-modules not containing $\mathbb C\Gamma$ and $W$ runs over $\hat \Gamma$,
the set of characters of $\Gamma$.

The advantage of working with $B^\tau$(rather than $O^\tau$) is that $B^\tau$ is
a ``one-relator'' algebra: it has a presentation as a quotient of a quasi-free
algebra $R$ by a two-sided ideal generated by a single element (see \eqref{rel1}).
Following \cite{BC}, we can think of this presentation as a differential
graded resolution of $B^\tau$. To be precise, let $ \boldsymbol B $ denote the graded associative 
algebra $\, I \oplus R\,$ having two nonzero components: the algebra 
$ R = \mathbb C \langle x,y\rangle \ast \Gamma $ in degree zero and its
(two-sided) ideal $ I := R\nu R $ in degree $ -1\,$. The differential on $ \boldsymbol
B $ is defined by the natural inclusion $\, d : I \hookrightarrow R
\,$ (so that $\,d\nu = xy-yx-\tau \in R \,$ and $\,da \equiv 0\,$
for all $\,a \in R\,$). Now there exists a canonical quasi-somorphism of $DG$ algebras
given by the projection $\eta : \boldsymbol B \to B^\tau$. This map yields the restriction functor
$\eta_{*} : \Mod(B^\tau) \to \DGMod(\boldsymbol B)$, which is an exact embedding.
It is well-known (see \cite{K}) that at the level of derived categories this functor
induces an equivalence of triangulated categories $D(\Mod(B^\tau)) \to
D(\DGMod(\boldsymbol B))$. 

Now, let $M$ be a projective $B^\tau$-module representing a class in $\mathcal R$.
We will associate to $M$ an object $\L$ in $\DGMod(\boldsymbol B)$ together with a 
quasi-isomorphism $M\to \L$, which we will call ${\it DG\!-\!model}$ of $M$. The 
DG-models are characterized by simple axioms (see Definition $2$ in Section $3.1$ below),
which determine $\L$ for each $M$ uniquely up to isomorphism 
(see Theorem \ref{uniq} in Section $5.2$). Thus our first result is following
\begin{theorem}
Let $\mathcal M$ be the set of isomorphism classes of DG-models in $\DGMod(\boldsymbol B)$. 
Then taking cohomology $\L\to H^{\bullet}(\L)$ induces a bijection 
$\omega_1: \mathcal M \iso \mathcal R$.
\end{theorem}
Next, in Section $3.2$, we will show that each DG-model determines a point in the union 
of Nakajima varieties:
\begin{equation}
\la{NVU}
\mathfrak M^{\tau} =  \bigsqcup_{V,W} \mathfrak M^{\tau}(V, W) \, ,
\end{equation}
where $V$ runs over the set of isomorphism classes of all finite dimensional 
$\Gamma$-modules and $W$ over the set of one-dimensional ones. Conversely, there is an 
explicit construction assigning to each point in $\mathfrak M^\tau$ a DG-model in 
$\DGMod(\boldsymbol B)$ (see Section $3.3$). In this way we will establish
\begin{theorem}
\la{th4}
There is a natural bijection $\omega_2 : \mathfrak M^\tau \to \mathcal M$.
\end{theorem}
Combining Theorems $2$ and $3$ together we arrive at the following result
(originally due to Baranovsky, Ginzburg and Kuznetsov, \cite{BGK}):
\begin{theorem}
\la{bgk}
The isomorphism classes of projective $B^\tau$-modules of rank(=dimension) $1$ are
in $1\!-\!1$ correspondence with points of the Nakajima varieties $\mathfrak M^\tau$.
\end{theorem}

If compared with \cite{BGK}, our proof of Theorem \ref{th4} has two main advantages.
First, we construct the bijection $ \Omega : \mathfrak M^\tau \longrightarrow \mathcal M$
as the composition of two maps $\Omega = \omega_2\circ\omega_1$, each of which 
is easy to describe. As a result, we give a completely explicit description
of rank $1$, projective $B^\tau$-modules (and thence, the right ideals of $O^\tau$). 

To be precise, let $\{ W_0, W_1, ..., W_{m-1}\}$ be the complete set of irreducible 
representations of $\Gamma = \mathbb Z_m$ such that $W_n \cong \epsilon^n$, and 
let $\{ e_0, e_1, ..., e_{m-1}\}$ be the corresponding idempotents in $\mathbb C\Gamma \subset B^\tau$. 
Writing $\mathfrak M^\tau =\bigsqcup_{n=0}^{m-1} \mathfrak M^\tau(W_n)$, 
we denote by $\Omega_n$ the restriction of $\Omega :\mathfrak M^\tau \to \mathcal R$ 
to the $n$-th stratum $\mathfrak M^\tau(W_n)$. Then we have the following theorem 
which extends the main result of \cite{BC}.
\begin{theorem}
\la{desideal}
The map $\Omega_n : \mathfrak M^\tau(W_n) \to \mathcal R$ sends a point of
$\mathfrak M^\tau(W_n)$ represented by a quadruple $(\X, \Y, \i, \j)$ to the class
of the fractional ideal of $B^\tau$:
$$M=e_n \det(\Y-yI) B^\tau + e_n\kappa \det(\X-xI) B^\tau \, ,$$
where $\kappa$ is the following element
$$ \kappa = 1- \j(\Y-yI)^{-1}(\X-xI)^{-1}\i(e_n)$$
in the classical ring of quotients of $B^\tau$. 
\end{theorem}
One of the interesting features of the Calogero-Moser correspondence
in the case of the Weyl algebra is its equivariance with respect to the action of the 
automorphism group $Aut(A_1)$. Our approach allows us to extend this result to the case of 
noncommutative Kleinian  singularities as follows. Let $G$ be the group of $\Gamma$-equivariant
automorphisms of the algebra $R=\mathbb C \langle x, y \rangle \ast \Gamma$, preserving the
the element $ xy-yx \in R$. For each $\tau \in \mathbb C \Gamma$, the canonical projection 
$R \to B^\tau$ yields a group homomorphism $G \to Aut_{\Gamma}(B^\tau)$ , and thus we have 
an action of $G$ on the space $\mathcal R$ (induced by twisting the right $B^\tau$-module
structure by automorphisms of $B^\tau$). On the other hand, there is a natural action of $G$
on the Nakajima varieties $\mathfrak M^\tau(V,W)$.
Finally, we observe that each $\sigma \in G$ extends naturally to an automorphism of the 
DG-algebra $\boldsymbol B$ and thus defines an auto-equivalence $\sigma_{*}$ on the category
$\DGMod(\boldsymbol B)$. It is easy to see that our axiomatics of DG-models is invariant under
such auto-equivalences and hence we have the induced action of $G$ on $\mathcal M$.
Now, the two bijections $\omega_1 : \mathcal M \to \mathcal R$ and $\omega_2 : \mathfrak M^\tau
\to \mathcal M$ obviously commute with the actions of $G$ defined above. Thus, we have the following
\begin{theorem}
\la{th5}
The map $\Omega : \mathfrak M^{\tau} \to \mathcal R$ is $G$-equivariant .
\end{theorem}
We remark that if $\Gamma =\{ e \}$ the group $G$ is isomorphic to $Aut_{\mathbb C}(A_1)$
(by a result of Makar-Limanov, \cite{ML}) and in this case our Theorem \ref{th5}
becomes one of the main results of Berest-Wilson.
In general, comparing our results with \cite{BW1} suggests that the (non-empty) subvarieties 
$ {\mathfrak M}^{\tau}(V, W)$ in \eqref{NVU} are precisely the orbits of the given action of 
$ G $ on $ \mathfrak M^{\tau}$. We will verify this conjecture in our subsequent paper.

\section{ K-theory}
The purpose of this section is to give a $K$-theoretic classification 
of ideals of $O^\tau$, that is the classification of ideals of $O^\tau$ up to stable isomorphism. 
Let us remind that we denote by $\mathcal R'$ the set 
of isomorphism classes of ideals of $O^\tau$ and by $\mathcal R$ the set
of isomorphism classes of $B^\tau$-submodules of $eB^\tau$. By Lemma $1$, 
these sets are in natural bijection. We will construct a map $\gamma: \mathcal R \to
 K_0(\Gamma) \times \hat \Gamma$ such that for any two isomorphism 
classes $(M_1)$ and $(M_2) \in \mathcal R$, the modules $M_1$ and 
$M_2$ are stably isomorphic if and only if $\gamma ((M_1))=\gamma ((M_2))$.

First, we would like to make some remarks about the Grothendieck groups $ K_0(\Gamma)$, 
$K_0(B^\tau)$ and $ K_0(O^\tau)$. We denote by $[\,\cdot \,]$, $[\,\cdot \,]_{K_0}$ 
and $[\,\cdot \,]_{K_0(O^\tau)}$ the stable isomorphism classes in the respective 
$K$-groups. By a well-known theorem of Quillen, the functor $ P \mapsto 
P\otimes_{\mathbb C\Gamma} B^\tau$ gives an isomorphism of groups 
$K_0(\Gamma) \cong K_0(B^\tau)$, and since the set $\{[W_n]\}^{m-1}_{n=0}$ 
generates $K_0(\Gamma)$, $\{ [e_n B^\tau]_{K_0}\}$ are the generators of $K_0(B^\tau)$. Furthermore, 
since $B^\tau$ and $O^\tau$ are Morita equivalent algebras, the corresponding equivalence 
functor induces another isomorphism $K_0(B^\tau) \cong K_0(O^\tau)$.
We will use these isomorphisms to identify $K_0(B^\tau)$ and $K_0(O^\tau)$ 
with $K_0(\Gamma)$. Now, the map assigning to a finite-dimensional module of $\Gamma$ 
its dimension extends to a group homomorphism $ \dim : K_0(\Gamma) \to \mathbb Z$
and have the following result:
\begin{proposition}
\la{neb}
Under the above identification of $K_0(\Gamma)$ and $K_0(O^\tau)$, the dimension function
coincides with the rank function on projective modules $\mathtt{rk}: K_0(O^\tau) \to \mathbb Z$.
\end{proposition}
\begin{proof}
Let $\widehat{\mathtt{dim}}: K_0(O^\tau)\to \mathbb Z$ be the composition of the isomorphism 
$K_0(O^\tau)\cong K_0(\Gamma)$ with $\mathtt{dim}: K_0(\Gamma) \to \mathbb Z$. We need 
to show that $\widehat{\mathtt{dim}}=\mathtt{rk}$. It suffices to check this on generators of $K_0(O^\tau)$,
say  $\{[e_n B^\tau e]_{K_0(O^\tau)}\}$. By definition of the dimension 
function we have $\widehat{\mathtt{dim}}([e_n B^\tau e])=1$. On the other hand, each of
the $e_n B^\tau e$ can be embbeded into $O^\tau$ as an ideal and therefore
$\mathtt{rk}(e_n B^\tau e)=1$.  
\end{proof}

Let us mention the following important result due to Baranovsky, Ginzburg and Kuznetsov
 (see  \cite{BGK}, Proposition $1.3.11$)
\begin{lemma}
\la{P1}
Let $ P\,\in \, K_0(\Gamma) $ be such that $ \mathtt{dim}(P)=1$. Then there exist $\Gamma$-modules
$W$ and $V$ (uniquely determined by $P$ up to isomorphism ), such that  we have 
$ P = [W] + [V] \cdot ([L]-2[W_0]) $ in $ K_0(\Gamma)$  and moreover,
$ \mathtt{dim}(W)=1$ and $ V $ does not contain the regular representation as a submodule.
\end{lemma}

Thus, with above identification of $K_0(B^\tau)$ and $K_0(\Gamma)$, we can define a map 
$\gamma: \mathcal R \to K_0(\Gamma)\times\hat \Gamma$ by:
$(M) \mapsto (V,W)$, where $(V,W)$ is the pair of $\Gamma$-modules from Lemma\! \ref{P1}.  
 Now, restating this Lemma  in terms of $B^\tau$-modules gives a classification modules in 
$\mathcal R$(or equivalently, in $\mathcal R$) up to stable isomomorphism.
\begin{theorem}
\la{Kth}
For any two isomorphism classes $(M_1), (M_2) \in \mathcal R$, we have 
$[M_1]_{K_0}=[M_2]_{K_0}$ if and only if $\gamma ((M_1))=\gamma ((M_2))$.  
\end{theorem}

Now, we will give a construction of the map $\gamma$ by
showing how to explicitly determine the $\Gamma$-modules $V$ and $W$ for
given class $(M)\in \mathcal R$.

Filter $B^\tau$ by assigning degree $1$  to the generators $x$ and $y$ and degree $0$ to 
all elements of $\Gamma$. Let us denote by $\bar{B^\tau}$ the
associated graded algebra and let $\bar M$ be the associated graded
module of a module $M \in \mathtt{mod}(B^\tau)$ equipped with a good
filtration. Each ideal $M$ of $B^\tau$ can be equipped with the
 induced filtration (which is good as $\bar B^\tau$ is Noetherian).
\begin{proposition}
\la{psi} For any isomorphism class $(M)$ in $\mathcal R$, there is
unique $n\in \{ 0,1,..., m-1\}$ such that $\bar M \hookrightarrow e_n
\bar{B}^\tau$ and  $\mathtt{dim}_{\mathbb C}(e_n \bar B^\tau /\bar
M)<\infty$.
\end{proposition}

The quotient $ e_n \bar B^\tau /\bar M$ can be viewed  as a
(finite-dimensional) $\Gamma$-module via the canonical inclusion
$\mathbb C\Gamma \to \bar B^\tau$.
\begin{lemma}
\la{lsi} 
Let $(M_1), (M_2)\in \mathcal R$ be such that $\bar M_1
\hookrightarrow e_n \bar{B}^\tau$ and $\bar M_2 \hookrightarrow e_k
\bar{B}^\tau$ with finite dimensional quotients, for some $n\in \{ 0,1,..., m-1\}$, then
$[M_1]_{K_0}=[M_2]_{K_0}$ if and only if $n=k$ and $[e_n \bar{B}^\tau / \bar 
M_1]=
[e_k \bar{B}^\tau / \bar M_2]$.
\end{lemma}

This lemma allows us to give an explicit construction of the map
$\gamma$. Specifically, let $(M)\in R$ be such that  $\bar M
\hookrightarrow e_n \bar{B}^\tau$ and $\mathtt{dim}_{\mathbb C}(e_n
\bar B^\tau /\bar M)<\infty$, then we can assign to the class of $M$
the pair $([W_n], [e_n \bar B^\tau /\bar M]) \in K_0(\Gamma) \oplus
K_0(\Gamma)$. We will show that this map coincides with
$\gamma$.

Let $G_0(\bar B^\tau)$ be the Grothendieck group of finitely
generated modules over $\bar B^\tau$. Then it is well-known (see,
for example \cite{G} Corollary $\!1.3$), that the class of $\bar M$
in $G_0(\bar B^\tau)$ does not depend on a choice of good filtration
on $M$, thus defining a map $\psi : K_0(B^\tau) \to G_0(\bar
B^\tau), [M] \mapsto [\bar M]_{G_0}$. Now, since both $\mathbb C\Gamma$
and $\bar B^\tau \cong \mathbb C[x,y] \ast \Gamma$ are Noetherian
rings of finite global dimension, $\psi$ is an isomorphism of groups. 
Moreover, we have the following commutative diagram:
\begin{equation}
\la{maps}
\begin{diagram}[small, tight]
        &                 &K_0(\Gamma)  &
       &    \\
        &\ldTo^{\phi_1} &   &\rdTo^{\phi_2}   \\
K_0(B^\tau) &                 &\rTo^{\psi}   & & G_0(\bar B^\tau) \, ,\\
\end{diagram}
\end{equation}
where $\phi_1$ and $\phi_2$ are group isomorphisms induced via
canonical embeddings of $\mathbb C\Gamma$ in the algebras $B^\tau$ and $\bar B^\tau$
respectively. We have
\begin{lemma}
\la{lex}
 Let $L$ be the natural two-dimensional representation of
$\Gamma$ and let $V$ be a finite-dimensional module over $\bar
B^\tau$. Then there is following class equation in $K_0(\Gamma)$ :
\begin{equation}
\la{classequ}
 \phi_2 ^{-1}([V]_{G_0}) = [V](2[W_0] - [L]),
\end{equation}
\end{lemma}
\begin{proof}
We consider the following sequence of $ \bar B^{\tau}$-modules
\begin{equation}
\la{dfg}
\begin{diagram}[small, tight]
                              & (V \otimes \epsilon)
                              \otimes_{{\mathbb C}\Gamma}\bar B^{\tau}
                              &                        \\
0 \to V \otimes_{{\mathbb C}\Gamma}\bar B^{\tau}\stackrel{d_2}
{\longrightarrow}  &  \oplus & \stackrel{d_1}{\longrightarrow} V
\otimes_{{\mathbb C}\Gamma}\bar B^{\tau}  \stackrel{d_0}
{\longrightarrow} V \to 0                               \\
                                & (V \otimes \epsilon^{-1})
                                \otimes_{{\mathbb C}\Gamma}\bar B^{\tau}  &
             \\
\end{diagram}
\end{equation}
where the maps are given by
\begin{eqnarray}
&& d_0 (v_1 \otimes b_1) = v_1.b_1 \, , \nonumber\\*[1ex] 
&&  d_1(v_1 \otimes \epsilon \otimes b_1,\, v_2 \otimes \epsilon^{-1} \otimes
b_2) = (v_1.y \otimes b_1 - v_1\otimes y \cdot b_1) - (v_2.x \otimes b_2 - v_2 \otimes x \cdot b_2)
\nonumber\\*[1ex] 
&&
d_2 (v_1 \otimes b_1) = (v_1.x \otimes \epsilon \otimes b_1 - v_1\otimes \epsilon \otimes x \cdot b_1\, , \,
v_1.y \otimes \epsilon^{-1} \otimes b_1 - v_1\otimes \epsilon^{-1} \otimes y \cdot b_1) \nonumber
\end{eqnarray}
for $v_i \in V $ and $b_i \in \bar B^{\tau}$ ($ i=1,2$).
We claim that this sequence is exact. First, it is easy to see that
$ d_2 \circ d_1=d_1 \circ d_0=0$. Second, it is clear that $d_0$ is
surjective and that $ \mathtt{Ker}(d_0) = \mathtt{Im} (d_1) $. So we
only need  to prove that $d_2$ is injective and that
$\mathtt{Ker}(d_1) = \mathtt{Im}(d_2)$. 

Let $\{ v_1,...,v_n \} $ be a basis of the finite-dimensional space $V$ and let $\bar X =(X_{ij})$
and $\bar Y=(Y_{ij}) $ be matrices corresponding to the actions of $x$ and $y$ in this basis.
Now if $u=\sum_{i=1}^{n} v_i\otimes b_i \in \mathtt{Ker}(d_2)$, then
$ x b_i = \sum_{j=1}^{n} X_{ij} b_j $. Assumig that $u \not= 0$ we let
$b_{i_0}$ be the element of largest degree among $\{ b_1,...,b_n \}
\subset \bar B^{\tau}$ with respect to the above filtration. Then, $
\mathtt{deg}(x b_{i_0}) >  \mathtt{deg}(\sum_{j=1}^{n} X_{i_0j} b_j)$ which contradicts to
the above equality and therefore $u=0$, so $d_2$ is injective.

Now if $ (u,u') = (\sum_{i=1}^{n} v_i \otimes \epsilon \otimes b_i,\, \sum_{i=1}^{n}
v_i\otimes \epsilon^{-1}\otimes c_i) \in \mathtt{Ker(d_1)} $, then 
\begin{equation}
\la{kerd1} \sum_{i=1}^{n} v_i \otimes ( \sum_{j=1}^{n} Y_{ij} b_j -
y b_i ) = \sum_{i=1}^{n} v_i \otimes ( \sum_{j=1}^{n} X_{ij} c_j - x \, ,
c_i )
\end{equation}
and we have 
\begin{equation}
\la{discrete} \sum_{j=1}^{n} Y_{ij} b_j - y b_i = \sum_{j=1}^{n}
X_{ij} c_j - x c_i , \quad i=1, ..., n \, ,
\end{equation}
which we can simply write as follows: $ (\bar Y-yI) \boldsymbol b = (\bar X -xI) \boldsymbol c$,
 where  $\boldsymbol b$ and $\boldsymbol c$ are column vectors consisting
of $b_i$ and $c_i$, $i=1, ..., n,$ respectively.

To prove that $\mathtt{Ker}(d_1) = \mathtt{Im}(d_2)$ we must show that that there 
exists $ u'' = \sum_{i=1}^{n} v_i \otimes d_i $ such that $u = (\X - xI) u''$ and $u' = (\Y-yI) u''$.
This is equivalent to finding a column vector $\boldsymbol d$  consisting of $d_i$, $i=1, ..., n,$ 
such that $ \boldsymbol b = (\bar X -xI) \boldsymbol d$ and $\boldsymbol c =(\Y -yI) \boldsymbol d$. 
From the matrix equation of \eqref{discrete} we can derive that each of $b_i$ is divisible by
$\mathtt{det} (\X -xI)$ and each $c_i$ by $\mathtt{det} (\bar Y-yI)$. 
Now, if we choose $\boldsymbol d:=(\X - xI)^{-1} \boldsymbol b=(\Y -yI)^{-1} \boldsymbol c$, 
then it satisfies the required property.
This proves exactness of the sequence \eqref{dfg}.

Thus, from \eqref{dfg} we obtain the following class equation in
$G_0(\bar B^{\tau})$:
\begin{eqnarray}
\la{main} 
&& [V]_{G_0}=[V \otimes_{{\mathbb C}\Gamma}\bar
B^{\tau}]_{G_0}- [(V\otimes \epsilon)\otimes_{{\mathbb C}\Gamma}\bar
B^{\tau}]_{G_0} \\*[1ex]
&&\quad  -[(V \otimes \epsilon^{-1})
\otimes_{{\mathbb C}\Gamma}\bar B^{\tau}]_{G_0} + [V
\otimes_{{\mathbb C}\Gamma}\bar B^{\tau}]_{G_0} \nonumber
\end{eqnarray}
Now applying $\phi_2^{-1}$ to \eqref{main} we get the desired
identity.
\end{proof}

\begin{proof}\!\!\!(of Lemma \!\ref{lsi})
We recall that $\psi: K_0(B^\tau) \to G_0(\bar B^\tau)$ is a group isomorphism and therefore
 $[M_1]_{K_0}=[M_2]_{K_0}$ if and only if $[\bar M_1]_{G_0}=[\bar M_2]_{G_0}$.
The inclusions $\bar M_1 \hookrightarrow e_n \bar B^\tau$ and $\bar M_2 \hookrightarrow 
e_k \bar B^\tau$ yield the following identities in $G_0(\bar B^\tau)$:
$$ [\bar M_1]_{G_0} = [e_n \bar B^\tau]_{G_0} - [e_n \bar B^\tau / \bar M_1]_{G_0} \quad
\mbox{and} \quad 
 [\bar M_2]_{G_0} = [e_k \bar B^\tau]_{G_0} - [e_k \bar B^\tau / \bar M_k]_{G_0} .$$
Applying  to these identities the group isomorphism $\phi_2^{-1}: G_0(\bar B^\tau) \to K_0(\Gamma)$ 
and using \eqref{classequ} we obtain
\begin{eqnarray}
\la{ko} 
&& \phi_2^{-1}([\bar M_1]_{G_0}) = [W_n] + [e_n\bar B^\tau / \bar M_1]([L] - 2[W_0]) \nonumber \\*[1ex]
&& \phi_2^{-1}([\bar M_2]_{G_0}) = [W_k] + [e_k\bar B^\tau / \bar M_2]([L] - 2[W_0]) \, ,\nonumber
\end{eqnarray}
and now the statement easily follows from Lemma \!\ref{P1}.
\end{proof}
\section{ DG-models}
\subsection{Axioms}
Let us remind that we denote by $ \boldsymbol B $ the graded associative algebra 
$\, I\!\oplus R\,$ having two nonzero components: the quasi-free algebra $ R = \mathbb C
\langle x,y\rangle \ast \Gamma $ in degree zero and its(two-sided)
ideal $ I := R\nu R $ in degree $ -1\,$. The differential on $ \boldsymbol
B $ is defined by the natural inclusion $\, d : I \hookrightarrow R
\,$ (so that $\,d\nu = xy-yx-\tau \in R \,$ and $\,da \equiv 0\,$
for all $\,a \in R\,$). The canonical map $ f : R \to \boldsymbol B
$ yields the restriction functor $f_{*} : \DGMod(\boldsymbol B) \to
\Com(R)$. So any DG module may be viewed as a complex of $R$-modules
and, in particular, as a complex of $\mathbb C \Gamma$-modules (via
the inclusion of $\mathbb C \Gamma$ into $R$).

We also recall that $\mathcal R (V,W)$ is the set of isomorphism classes
of finitely generated, projective (right) modules $M$ over $B^\tau$ such that
$[M]_{K_0} = [W] + [V]([L] - 2[W_0])$ in $K_0(B^\tau) \cong K_0(\Gamma)$ and $\mathtt{dim} [M] =1$.
Now if $M \in \mathcal R (V,W)$ for some finite dimensional $\Gamma$-module $V$
and  $W=W_n \in \hat \Gamma$, then we introduce the 
following definition (see \cite{BC}).
\begin{definition}
\la{D1}
A {\it DG
-model}\, of $ M $ is a quasi-isomorphism $\, q: M \to \LB \,$ in $
\DGMod(\boldsymbol B) $,
where  $\, \LB = L^0 \oplus L^1 \, $ is a DG-module with two nonzero
components
(in degrees $ 0 $ and $1$) satisfying the conditions:
\begin{itemize}
\item {\it Finiteness:}
\begin{equation}
\la{2.1}
\mathtt{dim}_{\mathbb C} L^1 \, < \, \infty .
\end{equation}
\item {\it Existence of a cyclic vector:}
\begin{equation}
\la{2.2} \mbox{There exists a $\Gamma$ linear map } i : W \to L^0
\mbox{ such that }  i(W).R = L^0 .
\end{equation}
\item {\it `Rank one' condition:}
\begin{equation}
\la{2.3}
\LB.\nu \subseteq  \mathtt{Im}(i)\ ,
\end{equation}
where $ \LB.\nu $ denotes the action of $ \nu $ on $ \LB\,$
and $\mathtt{Im}(i) \,$ denotes the image of   $ i \,$ in $ L^0 $.
\end{itemize}
\end{definition}
\noindent
The following properties are almost immediate from
the above definition.

$1.$ Since $W$ is a one-dimensional $\Gamma$-module there is a
canonical inclusion $W\hookrightarrow\mathbb C\Gamma$ under which $W
=  \mathbb C e_n$ . Then condition \eqref{2.2} says that $L^0$ is a
cyclic $R$-module with cyclic vector $i(e_n)$ which we denote by $
i_n$.

$2.$ The differential on $ \LB $ is given by a surjective $R$-linear
map: $ d_{\LB}:  L^0 \to L^1 $. This follows from \eqref{2.1} and
the fact that the cohomology of a DG module over $\boldsymbol B $ is a
complex of $B^\tau$ modules and that $B^\tau$ does not have
finite-dimensional modules (\cite{CBH} Th.$0.4$). Composing
$d_{\LB}$ with $i$ one obtains the map $ \i : W \to L^1$. Again as
in $1$, we can conclude that $L^1$ is a cyclic $R$ module with
cyclic vector $\i(e_n)$ which we denote by $\i_n$.

$3.$ Since $\nu$ is a degree $-1$ element in $\boldsymbol B$ we have
$L^0.\nu = 0$. Thus condition \eqref{2.3} in the axiom is equivalent
to $L^1.\nu \subseteq  Im(i) $. Define a map $\j : L^1 \to W $ so
that $v.\nu = i(\j(v))$ .
By Schur's Lemma the map $i$ is an injective map and therefore
$\j$ is a well-defined map . Now since $\Gamma \subset SL_2(\mathbb C)$
we have $ g \nu = \nu g$ for
all $g \in \Gamma$ which implies that $\j$ is a $\Gamma$-linear map.
Composing $\j$ with $d_{\L}$ we obtain another $\Gamma$-linear map
$ j : L^0 \to W $.

The following results give a useful characterization of  DG models
in the case of $\tau = 0$.
\begin{proposition}
\la{P2}
Suppose that $ \boldsymbol B_0 = I_{0} \oplus R $ ,
where $ I_{0} = R \nu_{0} R$, be a DGA such that
$d\nu_{0} = xy-yx $. If $ \LB =L^0 \oplus L^1 \in \DGMod(\boldsymbol B_{0})
$ satisfies
\eqref{2.1}-\eqref{2.3} then $L^1.\nu_{0} = 0 $ on $\LB$.
\end{proposition}
\begin{proof}
If $ L^1 = 0 $, then there is nothing to prove. So we may assume $L^1
\not= 0\,$. Then $ d_{\LB}(i) \not= 0 \,$ for  map $ f : R \to L^1
\,$ , $a \mapsto d_{\LB}(i).a $, is surjective by \eqref{2.2}. Now,
using the notation \eqref{2.4111} -- \eqref{2.6} and arguing as in
Lemma~\ref{L1} below, we can compute $\, [\X, \Y] = \i\,\j \,$. On
the other hand, the set of vectors $ \{\, \bar{G}\Y^m\X^k(\i) \,\} $
spans $ L^1 $ and $\,dim_{\mathbb C}\, L^1 < \infty \,$. An elementary lemma
from linear algebra (see, e.g., \cite{N}, Lemma~2.9) forces then $
\j = 0 $.
\end{proof}
\begin{corollary}
\la{CP2} Let $\boldsymbol B_0$ and $\LB$ be as Proposition~\ref{P2} then
$\LB$ can be identified with complex of $\bar B^{\tau}$-modules.
\end{corollary}
\subsection{The Nakajima data}
\la{Nd}
Let $ \LB $ be an DG-module satisfying the
axioms (\ref{2.1}) -- (\ref{2.3}). Denote by $\, X, Y , G\,$ (resp., $\, \X, 
\Y,
\bar G $) the action
of the canonical generators of $R $ on $ L^0 $ (resp., $ L^1 $), i.~e.
\begin{equation}
\la{2.4111} 
X(u) :=  u.\,x \in \End_{\mathbb C}(L^0) \ , \quad \X(v) :=  v.\,x \in \End_{\mathbb C}(L^1) \ ,
\end{equation}
\begin{equation}
\la{2.4222}
  Y(u) :=  u.\,y \in \End_{\mathbb C}(L^0)\ , \quad  \Y(v) :=  v.\,y \in \End_{\mathbb C}(L^1) \ ,
\end{equation}
\begin{equation}
\la{2.4333}
G(u) :=u.\, g \in \End_{\mathbb C}(L^0) \ ,  \bar G(v) :=v.\, g \in \End_{\mathbb C}(L^1) \ ,
\end{equation}
One can easily check that these maps satisfy the following
conditions
\begin{equation}
\la{2.6}
\X \,  d_{\LB} =  d_{\LB} \, X \ ,\quad \Y \,  d_{\LB} =  d_{\LB} \, Y \, 
\quad
\bar G\,  d_{\LB} =  d_{\LB} \, G .
\end{equation}
The next lemma shows that the linear data $(\X,\Y,\i,\j)$ extracted
from a DG model satisfy conditions \eqref{ma} and \eqref{mg} and
hence corresponds to a point in Nakajima variety.
\begin{lemma}
\la{L1}
The data introduced above satisfy the equations:
\begin{equation}
\la{2.9}
X\,Y - Y \,X + \it{T} = i\,j \ , \quad  \X\,\Y - \Y\,\X + \bar{T} = \i\, \j\
.
\end{equation}
\begin{equation}
\la{2.10}
X\,G = \epsilon (g) G\,X  , \quad  Y\,G = \epsilon(g) G \,Y , \quad
\X\,\bar{G} = \epsilon (g) \bar{G}\,\X  , \quad  \Y\,\bar{G} = \epsilon (g)
\bar{G}  \,\Y
\end{equation}
\end{lemma}
\begin{proof}
In view of \eqref{2.6} and surjectivity
of $ d_{\LB} \,$, the second parts of the equations 
(\ref{2.9}) and (\ref{2.10}) can be derived from
the first ones, and the first of (\ref{2.9}) follows easily
from the Leibnitz rule:
\begin{eqnarray}
&& T(u) = u. \, \tau = u. (xy-yx-d\nu) = u. xy - u. yx - u. d\nu =
\nonumber\\*[1ex] &&\quad   (u. x). y - (u. y). x + d_{\LB}(u). \nu
= Y\,X(u)  - X\,Y(u) + \,i (\j\,d_{\LB}(u)) = \nonumber\\*[1ex] &&
\quad  Y\,X(u)  - X\,Y(u) + i (j(u))=  (Y\,X - X\,Y + i\,j)\,u \ .
\nonumber
\end{eqnarray}
for all $\, u \in L^0$. 

Now to prove \eqref{2.10} we notice that 
$$ u. gx = (u. g). x = X\,G(u) \, ,$$
on the other hand we have
$$ u. \epsilon(g) x g = \epsilon(g) (u. x). g = \epsilon(g) G\,X(u). $$
\end{proof}
\subsection{From the Nakajima data to DG-models}
Let $\, (\X, \Y, \i, \j) \, \in End(U, U)^{\oplus 2} \oplus\\ 
Hom_{\Gamma}(W, U) \oplus Hom_{\Gamma}(U, W) $, where
$W\cong W_n$, be a quadruple representing a point in the Nakajima
variety. As a $\Gamma$-module $U$ can be uniquely written as 
$U \cong V \oplus \mathbb C \Gamma^{\oplus k}$ for some nonnegative
 integer $k$ and a module $V$ which does not contain the regular representation. 
If we denote $U$  by $L^1$, then, due to the stability condition $(ii)$ in \eqref{NV}, it is
clear that $L^1$ is a cyclic module with cyclic vector $\i(e_n) = \i_n$.

Using Nakajima data, we can introduce a functional
$\lambda:\, R \to \mathbb C $ so that $ \j(\i_n.\,a) =\lambda (a)\, \i_n$.
\begin{proposition}
\la{funck}
The functional $\lambda$ is defined by its values on the elements of the form
$x^ky^l$, where $k,l \ge0$. Moreover, we have $\lambda (x^ky^l)=0$ 
for all $k$ and $l$ such that  $k\not\equiv l(mod \,m)$.
\end{proposition}
\begin{proof}
First, by the condition \eqref{mg} and the fact that $g\in \Gamma$ acts on $i_n$ 
by a constant, it is sufficient to define $\lambda$ on the elements of the free
algebra $\mathbb C\langle x,y \rangle $. Second, due to the condition \eqref{ma}
we can express $\lambda (a)$, for any $a \in \mathbb C\langle x,y \rangle$,
in terms of $\lambda (x^ky^l)$, $k,l\ge0$. This finishes the proof of the first statetement.

Now, by definition of $\lambda$ we have
$$ \lambda(x^k y^l)\i_n=\j(\i_n.\,x^ky^l)=\j(\Y^l\X^k(\i_n))\ ,$$
and hence, using \eqref{alt} and \eqref{mg} sufficiently many times  we get
$$ \epsilon^n(g)\,\j(\Y^l \X^k(\i_n))=\j(\bar G\,\Y^l\X^k(\i_n))=\epsilon^{n+l-k}(g)\,\j(\Y^l \X^k(\i_n)) . $$
Finally, comparing the last two expressions we obtain the desired identity.
\end{proof}

Now we form the following {\it right} ideal in $ R $ :
\begin{equation}
\la{id}
J := \sum_{a \in R} (a \, (xy-yx-\tau) + \lambda(a))\,R\ .
\end{equation}
Let $\, L^0 := W\otimes_{{\mathbb C}\Gamma} R/J\,$, then 
since $\dim\,W=1$ we have that $ L^0 $ is a cyclic module over $R$ with the 
generator $\,e_ {n}\otimes [\,1\,]_J \,$ and hence we can define a map
$i : W \to L^0$ by $ w \mapsto w\otimes   [\,1\,]_J \,$.

If we consider the  map $ W \otimes_{\mathbb C} R \to L^1$ , $ w \otimes
a \mapsto  \i(w).\,a $, then elements of the form  $ w.g \otimes a - w \otimes
ga $ are annihilated by this map for any $ w \in W $ and $ a \in R $. Therefore
this map factors through the canonical projection $ W \otimes_{\mathbb C} R
\to W \otimes_{{\mathbb C}\Gamma} R $ inducing a map $ f : W
\otimes_{{\mathbb C}\Gamma} R \to L^1$. 

Further, it is easy to see that
\begin{equation}
\la{idealJ}
\i_n.\left[a(xy-yx-\tau) + \lambda(a)\right] = 0 \ , \quad \forall a \in
R  \, ,
\end{equation}
which allows to factor $ f $ through yet another canonical projection
$W \otimes_{{\mathbb C}\Gamma} R \to W \otimes_{{\mathbb C}\Gamma} R/J $
producing a map from $ L^0 $ to $ L^1 $. We denote this map by $ d_{\LB} $.
Being composition of $\Gamma$-linear maps
$d_{\LB}$ is also  $\Gamma$-linear .

Thus, we have constructed a complex of cyclic $ R $-modules
$$
\LB :=[\,0 \to L^0 \stackrel{d_{\LB}}{\longrightarrow} L^1 \to 0 \,] \ ,
$$
with differential $ d_{\LB} $. We want to enhance this complex with
a DG-module structure over $\boldsymbol B$.
For this it is sufficient to define the action of $\nu$ on $L^1$ and 
we define it as follows:
$( \i_n.\,a).\, \nu = -e_n \otimes [\,\lambda(a)\,]_{J}$.
Due to \eqref{idealJ} this action is well-defined and it is also clear that
$ L^1.\nu \subseteq Im(i) $.

Summing up, starting with Nakajima data $(\X, \Y, \i, \j) $
we have constructed a DG-module $\LB$ that
satisfies  all the axioms of Definition $1$.

Finally we have to show that $ \LB$ represents a rank $1$ projective
module over $B^{\tau}$ of an appropriate class in $K_0$.
\begin{lemma}
\la{L2}
Let $ \LB $ be a DG-module over $\boldsymbol B$ constructed above. Then its 
cohomology $H^0(\LB)$ is a finitely generated, projective module over $B^\tau$ such that 
$[H^0(\LB)]_{K_0} = [W_n]+[V]([L] - 2[W_0])$  in $K_0( B^{\tau} ) =K(\Gamma) $,
and consequently a representative of some class in $\mathcal R(V,W_n)$.
\end{lemma}
\begin{proof}
Let us fix some standard filtration on $ R $, say $R_{k} = \mbox{span} 
\{ x^p y^q g : p+q \leq k, g \in \Gamma \} $,
and put the induced filtration on $I$, so that $\grd(R) \cong R $ and $\grd(I) \cong
I_{0}$.
We then can filter the complex $ \LB \,$ as follows:
$\,L^0_{k} := i_n.R_{k} $ and $\, L^1_{k} := \i_n.R_{k} $.
Using \eqref{2.9} it is easy to see that the given DG-structure on $ \LB $
descends to the associated graded complex $\,
\grd(\LB) := \oplus_{n\geq 0} \LB_k/\LB_{k-1} \,$ making it
into a DG-module over $ I_{0} \oplus R$.
This module satisfies the same axioms \eqref{2.1}--\eqref{2.3}
as $ \LB $, and hence by Corollary~\ref{CP2}, we have the following
short exact sequence of $ \bar B^{\tau}$-modules
\begin{equation}
\la{first}
0 \to \overline  {H^0(\LB)} \to \bar L^0 \to \bar L^1 \to 0 ,
\end{equation}
where $ \bar L^0 := \grd(L^0)$ and $ \bar L^1 := \grd(L^1)$.
In particular, we have an isomorphism of $\bar B^\tau$-modules
 $\bar L^0 \cong W_n \otimes_{{\mathbb C}\Gamma}\bar B^{\tau} \,$. 
Passing from $\Mod(\bar B^{\tau})$ to  $\Mod(e\bar B^{\tau}e)$ via Morita
equivalence we see that $ \overline  {H^0(\LB)}e $ is a
submodule of $ W_n\otimes_{{\mathbb C}\Gamma}\bar B^{\tau}e$. 
The module $W_n \otimes_{{\mathbb C}\Gamma}\bar B^{\tau}e 
\cong e_n \bar B^\tau e$ can be identified with an ideal of  $e\bar B^{\tau}e$ 
and so can be $\overline { H^0(\LB)}e$. Thus $\overline {H^0(\LB)}e$ is a f.g., 
rank $1$, torsion-free module over $e\bar B^{\tau}e$. 
By standard filtration arguments all the above properties lift to 
$ H^0(\LB)e \,$ viewed as a module over the algebra 
$ O^\tau = eB^{\tau}e $. Now by Theorem $0.4$ of \cite{CBH},
 the  $gldim(O^\tau)=1$ and therefore $H^0(\LB)e \,$ is projective.
 Again, in view of the Morita equivalence between $ O^\tau$ and
 $B^{\tau} $, we conclude that $ H^0(\LB)$ is a projective $B^\tau$-module.

Now we need to show that $\phi^{-1}_1([H^0(\LB)]_{K_0}) = [W_n] + [V]([L] - 2[W_0])$ 
 which is equivalent, by \eqref{maps}, to showing that
$\phi_2^{-1}([\overline  {H^0(\LB)}]_{G_0})=[W_n] + [V]([L] - 2[W_0])$. From \eqref{first} we have
$[\overline {H^0(\LB)}]_{G_0}=[ \bar L^0]_{G_0} - [\bar L^1]_{G_0}$. Since 
$\bar L^0 \cong W_n\otimes_{{\mathbb C}\Gamma}\bar B^{\tau}$, we get
that  $\phi_2^{-1}([\bar L^0]_{G_0})=[W_n]$. Next we know that $\bar L^1$ is a finite-dimensional module 
over $\bar B^\tau$ isomorphic to $V \oplus \mathbb C
\Gamma^{\oplus k}$, and therefore by Lemma~\ref{lex} we obtain 
$\phi_2^{-1}([\bar L^1]_{G_0})=[V](2[W_0]-[L])$.
\end{proof}
\section{DG-models and Injective Resolutions}
\la{MRes} In this section we show how to construct some explicit
representatives of (the isomorphism class of) a module $ M$, such that 
 $(M)\in \mathcal R(V,W_n) $, from its $DG$-model $\,
M\stackrel{r}{\longrightarrow} \LB \,$. The key idea is to relate $
\LB $ to a minimal injective resolution of $ M \,$ (see \cite{BC}).

Let $\, \varepsilon: M \to \E \,$ be a minimal injective
resolution of $ M \,$ in $ \Mod(B^{\tau})\,$. 
Since global dimension of $B^\tau$ is one the resolution $\E$ has length one,
i.~e. $\, \E = [\,0 \to E^0 \stackrel{\mu_1}{\longrightarrow}  E^1
\to 0\,]\,$, and is determined (by $ M $) uniquely up to isomorphism
in $ \Com(B^{\tau})\,$. Recall that $\DGMod_{\infty}(\boldsymbol B)$ denotes the category
of DG modules over $\boldsymbol B$ with morphisms given by $\, \A$ homomorphisms.
Then, when regarded as an object in $\DGMod_{\infty}(\boldsymbol B)$, $\, \E \,$
is in the same quasi-isomorphism class as $ \LB \,$. 
It is natural to find a quasi-isomorphism that `embeds'
$ \LB $ into $ \E\,$. By Lemma \ref{L3} (see Appendix below) any 
$\, \A$-quasi-isomorphism between such modules is determined 
by two components $\, f = (f_1, f_2) $ where $ f_1 : \LB \to \E \,$ and
$ f_2 : L^1 \otimes  R \to E^ 0 $ .

\begin{theorem}
\la{T2} Let $\,r: M \to \L $ be a DG-model of $ M \,$, and let $\,
\varepsilon: M \to \E \,$ be a minimal injective resolution. Then
there is a unique $\, \A$-quasi-isomorphism $f_x : \L \to \E$ such
that   $(f_x)_1 \circ r =\varepsilon $ and
\begin{equation}
\la{3.16}
(f_x)_2\,(v, x) =0 \quad \mbox{and}
\quad (f_x)_2\,(v, g) = 0 \quad  \forall \, v \in L^1
\ ,\quad  \forall \, g \in \Gamma .
\end{equation}
\end{theorem}
\begin{remark}
First, a similar result can be stated if we replace $x$ by $y$. We will denote
the corresponding quasi-isomorphism by $f_y : \L \to \E$.
Second, the last equation of \eqref{3.16} implies that $f_2$ induces 
(and is determined by) the map $L^1 \otimes_{\mathbb C\Gamma} R \to E^0$
which we also denote by $f_2$. 
\end{remark}
\vspace{1ex}
The following lemma is essential for the proof of Theorem \ref{L2}.
\begin{lemma}
\la{P41} $ E^0 $ is a torsion-free module over  $ {\mathbb C}[x]\,$.
\end{lemma}
\begin{proof}
First of all , since $M$ is an ideal of $B^{\tau}$ we obtain that it is $ {\mathbb
C}[x]^{\Gamma}\,$-torsion free. Let $ n \in E^0 $ be a torsion
element then there is $q \in {\mathbb C}[x]^{\Gamma}\,$ such that
$q\neq 0$ and $ nq =0 $. Since $ E^0 $ is the injective envelope of $M$ we can find
nonzero $ b\in B^{\tau}$ and $m\in M$ such that $ m=nb $. Now the
elements of  $S= {\mathbb C}[x]^{\Gamma} \setminus \{0\} $ acts ad-nilpotently on
$B^{\tau}$ which implies that $S$ is an Ore set and hence there
are elements $t \in S$ and $c\in B^{\tau}$ such that $ bt = qc$ . 
Multiplying expression $ m=nb $ by $t$ we get :
$$ mt =nbt =nqc=0$$
which contradicts that $M$ is $ {\mathbb C}[x]^{\Gamma}\,$
torsion-free. This proves that $ E^0 $ is a torsion-free module over
$ {\mathbb C}[x]^\Gamma\,$. Now $\mathbb C[x]$ is a finite
integral extension of $ {\mathbb C}[x]^\Gamma\,$. 
Hence, for any nonzero $u \in \mathbb C[x]$ there exists a
minimal monic polynomial $f(v)=v^l+a_{l-1}(x)v^{l-1}+...+a_1(x)v+a_0(x)$ 
with coefficients in $ {\mathbb C}[x]^\Gamma\,$ such that $f(u)=0$
 and therefore we have
$$ u(u^{l-1}+a_{l-1}(x)u^{l-2}+...+a_1(x)) = -a_{0}(x).$$
If we had nonzero $n\in E^0$ such that $nu=0$ this would imply
$na_0(x)=0$ which contradicts that $E^0$ is torsion-free over $
{\mathbb C}[x]^\Gamma\,$.
\end{proof}

\begin{proof}(of $\mathbf {Theorem \, \ref{T2} }$)

First we observe that since there is a canonical inclusion of ${\mathbb
C}[x]\ast\Gamma$ into $R$ the complex $\L$ can be considered 
as a complex over ${\mathbb C}[x]\ast\Gamma$.
Now since $B^{\tau}$ is projective over ${\mathbb C}[x]\ast\Gamma$
(in fact, it is a free module $B^\tau = \oplus_{k=0}^{\infty} 
y^k {\mathbb C}[x]\ast\Gamma$), 
the complex $\E$ consists of ${\mathbb C}[x]\ast\Gamma$ injective 
modules. Hence, $ \varepsilon : M \to \E $ 
extends to a ${\mathbb C}[x]\ast\Gamma$-linear
morphism $\, f_1: \L \to \E \,$  such that the diagram
\begin{equation}
\la{D4}
\begin{diagram}[small, tight]
0 & \rTo &  M   & \rTo^{r}           & L^0         & \rTo^{d_{\L}} & L^1 & \rTo
   & 0 \\
  &      &  \dEq &                    & \dTo^{f_1}&            & \dTo^{\bar
f_1} &   \\
0 & \rTo &  M    & \rTo^{\varepsilon} & E^0         & \rTo^{\mu_1} & E^1&
\rTo   & 0 \\
\end{diagram}
\end{equation}
\noindent
commutes in $ \Com({\mathbb C}[x]\ast\Gamma)\,$. We claim that such an
extension is unique.
Indeed, if $\, f_{1}': L^0 \to E^0 \,$ is another map in $ \Mod({\mathbb
C}[x]\ast\Gamma) $
satisfying $\, f_1 \circ r = f_1' \circ r  = \varepsilon \,$,
then $\,  f_1' - f_1  \equiv 0 \,$ on $ \Ker(d_{\L})\,$ by exactness of the
first row of (\ref{D4}). So the difference
$\, \Delta := f_1' - f_1 \,$ induces a ${\mathbb C}[x]\ast\Gamma$-linear and
hence
${\mathbb C}[x]$-linear map $\, \bar{\Delta}: L^1 \to E^0 \,$. Since
$\, \dim_{\mathbb C}\,L^1 < \infty \,$,
$\, L^1 $ is torsion over $ {\mathbb C}[x]\,$, while $ E^0 $ is torsion-free
by the lemma above .
Hence, $\, \bar{\Delta} = 0 \,$ and therefore $\, f_1' = f_1 \,$.
This implies, of course, that $\, f_1'  = f_1 \,$ as morphisms in
$\Com({\mathbb C}[x]\ast\Gamma)$.
Now using part $b$ of Lemma \ref{L3} we can derive
that the map $\, f_1: \L \to \E \,$ extends to a
unique quasi-isomorphism of $\A$-modules over $ \boldsymbol B$ (see \cite{BC}
Lemma 6) .

Now the map $f$ being
${\mathbb C}[x]\ast\Gamma$-linear means that
\begin{equation}
f_2 (v, x) =0 \quad \mbox{and} \quad f_2 (v, g)  =0
\ ,\quad  \forall\, v \in L^1  ,\quad  \forall\, g \in \Gamma ,  
\nonumber
\end{equation}
which is exactly condition  \eqref{3.16}.
\end{proof}

To find the image of $M$ in $E^0$ we need to give an explicit
construction of $\,f_1 $. The formula \eqref{3.6} in Appendix B
relates $f_1$ and $f_2$. Substituting $c = \nu$ in this formula we
get
$$f_1(v.\nu) - f_1(v).\nu = -f_2(v, d\nu) \, .$$
Then as $\E$ is a complex over $B^{\tau}$ the second term
on the left hand side vanishes.
Now since $v.\nu= \j(v) i_n $ and $d\nu = xy -yx -\tau$, we obtain
the following equation on $f_2$ :
\begin{equation}
\la{3.19}
i_{x} \j(v) = -f_2(v, xy) + f_2(v, yx) -f_2 (v, \tau).
\end{equation}
where $\, i_x := f_1(i_n) \in E^0 \,$. Using \eqref{3.7} and
\eqref{3.16} we can rewrite this equation in the following form
\begin{equation}
\la{3.20}
f_2(v,y) \cdot x - f_{2}(\X(v), y)=  \j(v)\,i_x \ .
\end{equation}
Once this functional equation is solved one can recover $f_1$ from
\eqref{3.4} . The solution of \eqref{3.20} is given in Theorem \ref{T3} below.
To state this theorem we need the following important result.
\begin{lemma}
The ring $B^\tau$ has the classical (right) ring of quotients
$Q(B^\tau)$.
\end{lemma}
\begin{proof}
By Theorem $0.4$ of \cite{CBH}  the ring $O^\tau $ is simple and
therefore semiprime.
Being semiprime ring is a Morita invariant property so the ring $B^\tau$
is also semiprime. Now, as $B^\tau$ is a Noetherian, the existence
of $Q(B^\tau)$ is a consequence of Goldie'sTheorem (see \cite{St} pp.54-56)
\end{proof}
\begin{theorem}
\la{T3}
Let $ f_x $ be a $\A$-quasi-isomorphism defined in Theorem \ref{T2}, and let
$f_y$ be its counterpart obtained by interchanging $x$ and $y$ (see Remark 
following Theorem \ref{T2}). Then $\, f_x \,$ and $\, f_y \,$ are given explicitly by
\begin{eqnarray}
&&(f_x)_1\,(i_n.\, x^k y^m) =
i_x \cdot \left(\,x^k y^m + \Delta_{x}^{k m}(\i_n) \,\right)\ , \quad
(f_x)_2\,(v,\,x^k y^m) =  i_x\cdot \Delta_{x}^{k m}(v)\ , \la{3.18}\\*[1ex]
&&(f_y)_1\,(i_n.\,x^k y^m) = i_y \cdot \left(\,x^k y^m + \Delta_{y}^{k 
m}(\i_n)\,\right)\ , \quad
(f_y)_2\,(v,\,x^k y^m) =  i_y\cdot \Delta_{y}^{k m}(v) \ ,\la{3.181}
\end{eqnarray}
where $\, i_x := (f_x)_1\,(\,i_n\,) \,$ and $\, i_y :=
(f_y)_1\,(\,i_n\,)\,$ in $ E^0 $,  and
\begin{eqnarray}
&& \Delta_{x}^{k m}(v) := - \, \j(\X - xI)^{-1} (\Y-yI)^{-1} (\Y^m -
y^mI)\,\X^k v \ , \la{2.22x}\\*[1ex] && \Delta_{y}^{k m}(v) :=  \j(\Y
- yI)^{-1} (\X-xI)^{-1} (\X^k - x^kI)\,y^m v \ ,\la{2.22y}
\end{eqnarray}
where $I:=Id_{L^1}$. Moreover,
\begin{equation}
\la{iii} i_x .g =\epsilon^n(g)\, i_x , \quad i_y .g =\epsilon^n(g)\, i_y
,\quad  \forall\, g \in \Gamma \
\end{equation}
\begin{equation}
\la{ii}
i_x = i_y \cdot \kappa
\end{equation}
where $\kappa \in Q$ is given by the formula $\kappa=1 - \j\,(\Y - yI)^{-1} (\X-xI)^{-1}\i_n $  and 
satisfies the equation 
\begin{equation}
\la{com}
e_n \, \kappa \, (1 - e_n) = 0 \quad \mbox{ in } Q.
\end{equation}
\end{theorem}

\noindent
Let us give some comments on the theorem.

1.\ Since $i_n$ is a cyclic vector of one dimensional $\Gamma$-module
$W$ the elements $\{ i_n. \, x^k y^m \}$ form a basis of $L^0$. So it suffices 
to define the maps $(f_x)$ and $(f_y)$ only on these elements.

2.\ The formulas (\ref{2.22x}) and (\ref{2.22y}) define the maps $\,
\Delta_{x,y}^{k m}:\, L^1 \to Q(B^\tau) \,$ for $ m,k \geq 0\,$, which could
be written more accurately as follows
\begin{equation}
\Delta_{x}^{k m}(v) := - \det(\X - x\,\id)^{-1} (\j \otimes
1)\,[\,(\X - x \,\id)^{*}
\sum_{l=1}^{m}\, \Y^{m-l}\X^k(v) \otimes y^{l-1}\,] \ ,\nonumber
\end{equation}
where $\, (\X - x\,I)^{*} \in \End_{\mathbb C}(L^1) \otimes _{\mathbb C}  R \,$ 
denotes the classical adjoint of the matrix $\, \X - x\,I \,$ and 
$ \j \otimes1\,: L^1 \otimes _{\mathbb C} R \to R $ is defined by 
$\, v \otimes a \mapsto \j(v)\,a \,$.

3.\ The dot in the right hand sides of (\ref{3.18}) and \eqref{3.181}
denotes the (right) action of $ B^\tau $ on $ \E \,$. Even though
$\, \Delta_{x,y}^{k m}(v) \in Q(B^\tau) \,$, these formulas make sense
because both $E^0$ and $E^1$ are injective, and hence
{\it divisible} modules over $ B^\tau\,$.

\begin{proof}
The second formula of \eqref{3.18} and  can be checked simply by
plugging it in \eqref{3.20} for $v=\i_n$. Then using \eqref{3.4} we
derive the second formula of  \eqref{3.18}. Similarly considering
$\mathbb C[y] \ast \Gamma$ linear map $f_y$ we can obtain
\eqref{3.181}. The formulas in \eqref{iii} can be derived from the fact that both $f_x$
and $f_y$ are $\Gamma$-linear maps. 

Consider the polynomial $ p(x) =
\det(\X - xI)$ then by Hamilton-Cayley's theorem $ \i_n.\,p(x)=0. $
This implies that $ i_n.\,p(x) $ is in the image of $r$. Since $
f_x \circ r = \varepsilon = f_y \circ r$ we have $ f_x (i.\,p(x))
= f_y(i.\,p(x))$. Using \eqref{3.181} and \eqref{2.22y} we obtain
$$ i_x \cdot p(x) = i_y \cdot (1 - \j\,(\Y - yI)^{-1} (\X-xI)^{-1}\i_n)\, 
p(x)$$
and since $E^0$ is divisible module we derive formula \eqref{ii} simply dividing 
it by $p(x)$.

In order to prove \eqref{com} it sufficies to show that $e_n \kappa\cdot g =
\epsilon^n(g)\, e_n \kappa$ for all $g\in \Gamma$.
For this we expand $e_n \kappa$ into the formal series:
\begin{equation}
\la{ddf}
e_n \kappa = e_n - e_n \sum_{l, k \geq 0}\,
\j \left(\Y^l\X^k \,\i_n\right) \, y^{-l-1} x^{-k-1}=e_n-
e_n\sum_{l \equiv k(mod \, m)} \lambda_{kl}\,y^{-l-1}x^{-k-1}\, ,
\end{equation}
where $\lambda_{kl}=\lambda (x^k y^l)$ and the last equivality follows from Proposition \ref{funck}.
Now multiplying this series by $g$ we obtain
\begin{equation}
\la{mm}
e_n\kappa \cdot g = \epsilon^n(g)\, e_n \biggl (1 - \sum_{l \equiv k(mod \, m) }\,
\epsilon^{k-l}(g) \lambda_{kl} \, y^{-l-1} x^{-k-1} \biggr )=\epsilon^n(g)e_n\kappa \, .
\end{equation}
where the last equality follows from \eqref{ddf} and the fact that $\epsilon^{k-l}(g)=1$
for $l\equiv k (mod \, m)$. 
\end{proof}
\begin{corollary}
\la{C2}
Let $\, \L $ be an $DG$-envelope of $ M \in \mathcal R (V, W_n)$,
and let $\, (\X, \Y, \i, \j) \,$ be the Nakajima data associated with $ \L$. 
Then, $\, M \,$ is isomorphic to each of the
following (fractional) ideals
\begin{eqnarray}
&& M_x := e_{n} \det(\X - xI)\,B^{\tau} + e_{n} \mu
\,\det(\Y - yI)\,B^{\tau} \ , \la{2.Mx} \\*[1ex]
&& M_y := e_{n} \det(\Y - yI)\,B^{\tau}  +e_{n} \kappa
\,\det(\X - xI)\,B^{\tau} \ , \la{2.My}
\end{eqnarray}
where $\, \mu := 1 + \j\,(\X-xI)^{-1}(\Y - yI)^{-1} \i_n  $ is such that $e_n 
\kappa \cdot
e_n \mu = e_n \mu \cdot e_n \kappa = e_n$ and hence $M_y = e_n \kappa M_x$ .
\end{corollary}
\begin{proof}
Since $M$ is an ideal of $B^\tau$ and $E^0(M)$ is a divisible module
there is an inclusion $Q(B^\tau) \hookrightarrow E^0(M)$.
The key idea of the proof is to realize $M$ in its injective envelope 
$E^0(M)$
by investigating the image of $r(M)$  under the maps $f_x$ and $f_y$ .

Let $p(x) = \det (\X- xI)$ and $q(y) = \det (\Y-yI) $
then by Hamilton-Cayley's theorem  $ D = i_n.\,p(x)R + i_n.\,q(y)R$ 
is an $R$-submodule of $\mathtt{Ker}(d_{\LB})=\mathtt{Im}(r)$. 
It is easy to see that $D$ is a submodule of finite codimension in $L^0$
and hence $D$ is a finite codimension in $\mathtt{Im}(r)$.
Now, since $f$ is an injective map, $f_x(D)$ is a subspace of finite codimension in 
 $f_x (\mathtt{Im}(r))$. Further, it is clear that $f_x (\mathtt{Im}(r))=\varepsilon(M)$ 
is a $B^\tau$-module. 

If we show that $f_x (D)$ is also $B^\tau$-module, then, since for generic $\tau$ the algebra 
$B^\tau$ does not have finite-dimensional modules, we obtain $f_x (D)=f_x (\mathtt{Im}(r))$. 

By \eqref{3.18} and \eqref{2.22x} we have $f_x (i_n.\,p(x)R) = i_x \cdot 
p(x) B^\tau$.
Further, since $ f_x \circ r = \varepsilon = f_y \circ r$
we obtain $ f_x (i_n.\,q(y)R) = f_y(i_n.\,q(y)R) $ and
therefore by \eqref{3.181} and \eqref{2.22y} we have
$ f_y(i_n.\,q(y)R) = i_y \cdot q(y) B^\tau$.

Now, since
$$ [ \X - xI, \Y - yI] = [ \X, \Y] +  [x,y] I =
\i \circ \j -\bar T +\tau I \, ,  $$
it is easy to check that $e_n \mu \cdot e_n \kappa =
e_n \kappa \cdot e_n \mu =e_n $ and hence $i_y =i_x \mu $ .
Thus, we get
$$f_x (D)=f_x (i_n.\,p(x)R)+f_y(i_n.\,q(y)R)= i_x \det(\X - xI)\,B^{\tau} + i_x \kappa^{-1}
\,\det(\Y - yI)\,B^{\tau}. $$
By above arguments we obtain $M\cong \varepsilon(M)=f_x(\mathtt{Im}(r))=f_x(D)$ . 
To finish the proof we notice that there is a $B^\tau$-linear automorphism
of $E^0(M)$ sending $i_x$ to $e_n$.
\end{proof}
\section{Existence and Uniqueness}
\subsection{Distinguished Representatives}
\la{disrep}

In previous section, in Corollary \ref{C2}, for every $DG$ model $\L\in \mathcal M$
we have constructed two different realizations of  $H^0(\L)$ as fractional ideals of $B^\tau$. 
Our main goal in this section is to present analogous result for any $M\in \mathcal R$. 
This result will be essential for proving existence and uniqueness of $DG$ models.  

First, we notice that $S_1=\mathbb C[x] \setminus \{0\}$ is an Ore set in $B^\tau$.
Indeed we have already shown in
Proposition \ref{P41} that the set  $S=\mathbb C[x]^\Gamma \setminus \{0\}$
is an Ore set. Since $S$ is
an integral extension of $S_1$ then for any $u \in S_1$ we have
   $$ u(u^{k-1}+a_{k-1}(x)u^{k-2}+...+a_1(x)) = -a_{0}(x) $$
with $a_1(x), ..., a_{k-1}(x) \in \mathbb C[x]^\Gamma $
and $a_0(x) \in S$. So for any
$ b \in B^{\tau}$ there exist $c \in B^{\tau}$ and
$a \in \mathbb C[x]^{\Gamma}$ such that
$ab = ca_{0}(x)=[-c(u^{k-1}+a_{k-1}(x)u^{k-2}+...+a_1(x))] u $
which proves that  $S_1$
is an Ore set. Now let $B^{\tau}[S_{1}^{-1}]$ be the ring
of fractions of $B^{\tau}$  with
respect to $S_1$ then $M_x$ of Corollary \ref{C2} has the following
properties:
\begin{eqnarray}
\la{cond}
&& (1)\ M_x \subset e_nB^{\tau}[S_{1}^{-1}] \
\mbox{and}
\ M_x \cap e_n \mathbb C[x] \not= \{0\}\ , \nonumber \\
&& (2)\ \mbox{if} \ e_n(a_k(x) y^k + a_{k-1}(x)y^{k-1} + \ldots )
\,\in\, M_x \ \mbox{then} \
a_k(x) \in \mathbb C[x] \ ,\\
&& (3)\ M_x \ \mbox{contains elements of the form}\
e_n(y^k + a_{k-1}(x)y^{k-1} + \ldots) \ , \nonumber
\end{eqnarray}
We can also introduce the set $S_2=\mathbb C[y] \setminus \{0\}$ and show
that $M_y$ satisfies similar properties.
\begin{lemma}
Let $M$ be a representative of some class in $\mathcal R$, then there exists 
a fractional ideal $M_x$ of $B^\tau$ isomorphic to $M$ and satisfing conditions 
$(1)-(3)$ for some $n \in\{ 0,1,...,m-1\}$.
\end{lemma}
\begin{proof}
First of all, without loss of generality we may assume that $M$ is a submodule of
 $eB^\tau$ such that $M\cap e \, \mathbb C[x] \neq \{ 0\}$ (see \cite{BGK}, Lemma $6.4$).
Now let $\boldsymbol{w} =( w_1, w_2)$ be a pair of nonnegative real numbers
such that $\boldsymbol{w}= w_1+w_2 > 0 $ then
we introduce a natural increasing filtration on
$B^\tau$: $F_{\boldsymbol{w}}^{0} B^\tau = \mathbb C \Gamma ,
\, F_{\boldsymbol{w}}^{i} B^\tau =\{ x^k y^m g |\, \mathtt{deg}_{\boldsymbol{w}}(x^k
y^m g):=
kw_1 + mw_2 \leq i, g\in \Gamma \}$.
We can extend this filtration on $Q$ by the following formula:
$$ \mathtt{deg}_{\boldsymbol{w}} (ab^{-1}) :=
\mathtt{deg}_{\boldsymbol{w}}(a)-\mathtt{deg}_{\boldsymbol{w}}(b) $$
and $F_{\boldsymbol{w}}^{i} Q = \{ q\in Q |\, \mathtt{deg}_w(q) \leq i \}$.
With respect to this filtration we define
the associated graded algebra $\grd_{\boldsymbol{w}} B^\tau
=\oplus _{k=0}^{\infty} F_{w}^{k} B^\tau / F_{w}^{k-1} B^\tau \cong 
\mathbb C[\bar x,\bar y] \ast \Gamma$, where $\bar x := \grd (x)$
and $\bar y := \grd (y)$. We will now choose 
$\boldsymbol{w} = (0,1)$ and denote the associated graded
module of $M$ with respect to this filtration by $\grd_{y}(M)$.
Then we have 
$$ \grd_y (M) = \bigoplus_{k\ge 0} e_n\, I_k (\bar x) \, y^k , $$
where $I_0 \subseteq I_1 \subseteq I_2 \subseteq ... $ is an ascending
chain in $\mathbb C[\bar x]$ with $I_0(\bar x) \neq 0$. Since $\mathbb C[\bar x]$
is a PID each $I_k$ is cyclic and the sequence of ideals $\{ I_k\}$ stabilizes :
$ I_{n_0} = I_{n_0 + 1} = I_{n_0 + 2} = ... $ starting from some $n=n_0 \ge 0$ .
Denote by $p = p(\bar x)$ the principal generator of $I_{n_0} $ in 
$\mathbb C [\bar x]$. Now we claim that $p(\bar x) = x^j \tilde{p}(\bar x)$ for some 
$j \in \{0,1,...,m-1\}$ and $\tilde {p}(\bar x)$ is a $\Gamma$-invariant polynomial. 
It is clear that we can write $p(x)$ as follows: $p(x) = \sum_{j=0}^{m-1} x^j p_j(x)$, 
where each $p_j(x)$ is a $\Gamma $-invariant polynomial. 
Let $j_0$ be the minimal $j$ for which $p_j(x)\neq 0$ then there exists an 
element of $M$ of the form $b = e\, x^{j_0}  p_{j_0}(x) y^{n_0} +\{ \textit{terms of degree of } y
 \textit { less than } j_0\}$. Hence $\grd_y(b) = e\, \bar x^{j_0}  p_{j_0}(\bar x) \bar y^{n_0}$ which
implies our claim. Finally let  $M_x =p^{-1}(x) M$ then since $e x^{j_0} = x^{j_0} e_{m-j_0}$
we obtain $M_x \subset e_{m-j_0} B^\tau [S_1^{-1}]$
and $M_x$ satisfies conditions $(1)-(3)$. 
\end{proof}
\begin{corollary}
\la{4.1}
Let $M_{x}$ and $M_{x}'$ be two fractional ideals isomorphic
to $M$ and satisfying $(1)-(3)$ above. Let $q$  be an element of $Q$
such that $M_{x}' = qM_x$ then $ q\in \mathbb C e_n$ (and hence  $M_{x}' = 
M_x$) .
\end{corollary}
\begin{proof}
It is clear from $(2)$ of \eqref{cond} that
$\mbox{\rm gr}_{y}(M_x) \subset \grd_{y} B^\tau \cong e_n \mathbb C[\bar x,\bar y]$.
Moreover, due to conditions $(1)$ and $(3)$ this embedding is of finite
codimension. This in turn implies that  $\mbox{\rm gr}_{y} (q)
\in \mathbb C e_n$. Now since $ M_x \cap e_n \mathbb C[x] \not= \{0\}$
we have $q \in e_n \mathbb C(x)[y]$.
Combining these last two facts we conclude that $q\in \mathbb C e_n$.

\end{proof}

Reversing the roles of $x$ and $y$ , we obtain another
distinguished representative $M_y$.
The statement similar to Lemma~\ref{4.1} will also
be true for $M_y$. In Corollary~\ref{C2}
we have seen that there is an element $ \kappa \in Q$
such that $M_y = e_n \kappa M_x$. The following corollary
claims such element is unique.

\begin{corollary}
\la{Cor}
Let $M_x$ and $M_y$ be fractional ideals isomorphic to
$M$ and defined as above.
Let $q$  be an element of $Q$ such that $M_{y} = qM_x$ then
$ q$ is uniquely determined up to constant factor of $e_n$.
\end{corollary}
\begin{proof}
Suppose we have $q_1, q_2 \in Q$ such that $M_y = q_i M_x \,(i=1,2)$.
Since both $M_x$ and $M_y$ are submodules of $e_n Q(B^\tau)$  we
derive $q_1, q_2 \in e_n Q(B^\tau) e_n$.
Hence $M_x = qM_x$ where $q=q_2^{-1} q_1\in Q$. Now by the above
lemma $q\in \mathbb C e_n$.
\end{proof}

\subsection{Uniqueness}

In this section we will establish uniqueness of DG-models up to
isomorphism of DG-modules. First we remind the definition of linear
functional introduced earlier
$$   \lambda : R \to \mathbb C , \quad  a \mapsto
\lambda(a), \mbox{ where } \, \lambda (a)\, \i_n = \j(\i_n.\,a). $$
We recall from Section $3.3$ that $\lambda$ is completely deternined
by its special values: 
\begin{equation}
\la{numb}
\lambda_{kl} :=\lambda(x^k y^l) .
\end{equation}
\begin{theorem}
\la{uniq}
Let $\L$ and $ \tilde{\boldsymbol{L}}$ be two DG-models of $M$.
Then the following are
equivalent :

(a) $\L$ and $\tilde{\boldsymbol{L}}$ are isomorphic as DG-modules over  $
\boldsymbol B $,

(b)  $\L$ and $\tilde{\boldsymbol{L}}$ are $\A$-isomorphic,

(c) $\L$ and $\tilde{\boldsymbol{L}}$ are $\A$ quasi-isomorphic.

(d) $\L$ and $\tilde{\boldsymbol{L}}$ determine the same
functional $\lambda: R \to \mathbb C$
(i.e $\lambda = \tilde {\lambda} $).
\end{theorem}
\begin{proof}
The implications $\, (a) \Rightarrow (b) \Rightarrow (c) \,$ are obvious.
It suffices only to show that $ \, (c) \Rightarrow (d) \, $
and $ \, (d) \Rightarrow (a) \,$.

If $ \L $ satisfies \eqref{2.1}--\eqref{2.3} then, by Lemma~\ref{L2},
the cohomology $ H^0(\L) \in \mathcal R(V, W_n) $. By Corollary~\ref{C2},
$ H^0(\L) $ is then isomorphic to the fractional ideals $M_x$ and $M_y$
related by $\, M_y = e_n \kappa\,M_x \,$ (see \eqref{2.Mx}, \eqref{2.My}).
Expanding $e_n \kappa$ into the formal series as in \eqref{ddf} we have :
\begin{eqnarray}
\la{series}
&& e_n\,(1 - \j\,(\Y - y)^{-1} (\X-x)^{-1}\,\i) = e_n -e_n \sum_{l, k \geq
0}\,
\left(\,\j \,\Y^l\X^k \,\i\,\right) \, y^{-l-1} x^{-k-1}   \\*[1ex]
&& \ \ \ \ \ \ \ \ \ \ \ \ \ \ \ \ \ \ \  = e_n - e_n\, \sum_{l, k \geq 0}\,
\lambda_{kl} \, y^{-l-1} x^{-k-1} \nonumber
\end{eqnarray}
Now, $ e_n \kappa $ is determined uniquely,
up to a  constant factor of $e_n$,
by the isomorphism class of $ H^0(\L)$ (see Corollary \ref{Cor}).
Hence, if $ \L $ are $ \tilde{\boldsymbol{L}} $
are quasi-isomorphic $ \A$-modules,
we have $\,H^0(\tilde{\boldsymbol{L}}) \cong H^0(\L) \,$ and therefore
$\, e_n \tilde{\kappa} = c \cdot e_i\kappa \,$ for
some $ c \in \mathbb C\,e_n$. Comparing the
coefficients of \eqref{series} yields at once $\, c = e_n \,$ and
$\, \tilde{ \lambda}_{lk} =  \lambda_{lk} \,$ for all $ l,k \geq
0\,$.
Thus, we conclude $\,(c) \Rightarrow (d) \,$.

Now let $\lambda = \tilde{\lambda}$ then
$J = \tilde {J} $ and therefore the map $f_{1}^0 : L^0 \to \tilde {L}^0$ defined as
$ i_n.\,a \mapsto \tilde {i}_n.\,a $ is an isomorphism  of $R$-modules.
Since $\,H^0(\tilde{\boldsymbol{L}}) \cong H^0(\L) \,$
the induced map $f_{1}^1 :  L^1 \to \tilde {L}^1$
is also isomorphism. Finally since
$d_{\tilde {L}} \circ f_{1}^0 = f_{1}^1 \circ d_L$ the pair
$(f_{1}^0, f_{1}^1)$ produces the necessary DG-isomorphism
which proves implication $(d) \Rightarrow (a)$.
\end{proof}

\subsection{Existence}
\la{exist}
Let us start by stating the main result of the section.
\begin{theorem}
\la{existence} For each projective module $M\in \mathcal R (V, W_n)$
there is a DG model satisfying axioms of Definition $1$.
\end{theorem}
So we need to produce a right DG module which is a two-complex of vector 
spaces, quasi-isomorphic to $M$ and satisfies conditions \eqref{2.1} - \eqref{2.3}. 
Constructing such DG module $\L$  is equivalent to constructing a DG algebra 
homomorphism $ \Psi$ from $\boldsymbol B$ to the opposite of the DG algebra 
$\underline{Hom}_{\mathbb C}(\L,\L)$ which we denote by $\boldsymbol C$.
Let us remind that DG structure on $\boldsymbol C$. The multiplication is just
usual composition of endomorphisms and the differential is defined  by:
$$d_{\boldsymbol C}(f) := d_{\L} \circ f - (-1)^j f \circ d_{\L} ,
\, \mbox{ for } f \in C^j .$$
The algebra $\boldsymbol B$ has the generators $x, y, g(\in Z_m)$ in degree
zero and one generator $\nu$ in degree minus one such that $d_{\boldsymbol
B}(\nu) = xy - yx -\tau$. Since the algebra $\boldsymbol C$ consists of
three nonzero components $\boldsymbol C = C^{-1} \oplus C^0 \oplus C^1$
such that :
$$ C^{-1} = Hom_{\mathbb C}(L^1,L^0)$$
$$ C^0 =  End_{\mathbb C}(L^0)\oplus   End_{\mathbb C}(L^1)$$
$$ C^{1} = Hom_{\mathbb C}(L^0,L^1)$$
we need to map $x \mapsto (X, \X)$, $y \mapsto (Y, \Y)$ ,\, $g
\mapsto (G, \bar G)$ and $\nu \mapsto f$, where $X, Y, G \in
 End_{\mathbb C}(L^0)$,\quad $\X, \Y, \bar G \in  End_{\mathbb C}(L^1)$ 
and $f \in Hom_{\mathbb C}(L^1,L^0)$. Moreover,$\Psi$  must satisfy 
the following conditions $d_{\boldsymbol C}
(\Psi (x)) = d_{\boldsymbol C} ( \Psi (y)) = 0$ and $\Psi (d_{\boldsymbol B}\nu)=d_{\boldsymbol C}(\Psi(\nu))$ 
which are equivalent to the following system of equations:
\begin{equation}
\la{psi1}
d_{\L} \circ X = \X \circ d_{\L} \, , \quad d_{\L} \circ Y = \Y \circ d_{\L}
\end{equation}
\begin{equation}
\la{psi2}
X\,Y - Y \,X + T = f_0 \ , \quad  \X\,\Y - \Y\,\X + \bar T = f_1\
\end{equation}
where $(f_0, f_1)=d_{\boldsymbol C} (f)$. The rest of this section focuses
on the construction of such complex $\L$.

In Section~\ref{disrep} we have shown that ideals $M_x$ and $M_y$
defined in Corollary \ref{C2} are uniquely characterized by
properties $(1) - (3)$ of \eqref{cond}. Moreover,  by Corollary
\ref{Cor}, there is $e_n \kappa \in Q$ such that $M_y = e_n \kappa
M_x$ and $e_n \kappa$ is uniquely defined up to constant factor of $e_n$ . We
choose $e_n \kappa$ such that $\grd_y (e_n \kappa) = e_n$. Now even
though $M_x$ and $M_y$ are fractional ideals we can embed them into
$ e_n B^{\tau}$ by means of the following maps:
\begin{equation}
\la{rhox}
\rho_{x} : e_n  B^{\tau}[S_{1}^{-1} ] \to e_n B^{\tau} \, ,
\quad e_n b(x) y^m \mapsto e_n b(x)_{+} y^m
\end{equation}
\begin{equation}
\la{rhoy}
\rho_{y} :e_n  B^{\tau}[S_{2}^{-1} ] \to e_n B^{\tau} \, ,
\quad e_n b(y) x^m \mapsto e_n b(y)_{+} x^m
\end{equation}
where $"+"$ stands for taking polynomial part of corresponding rational function.
Let $r_x : M_x \to e_n B^\tau $ and $r_ y : M_y \to e_n B^\tau$ be restrictions of the 
above maps to $M_x$ and $M_y$ correspondingly and let $V_x = e_n B^\tau / r_x (M_x)$
and $V_y = e_n B^\tau / r_y (M_y)$. It is not difficult to see first that $r_x$
is $\mathbb C[y]\ast \Gamma$-linear and $r_y$ is $\mathbb C[x]\ast \Gamma$-linear maps
and second that both $V_x$ and $V_y$ are finite dimensional
$\Gamma$-modules. Now let's consider the following complexes of $\Gamma$-modules:
\begin{equation}
\L_x = : [\, 0 \to e_n B^\tau \to V_x \to 0 \,] \quad \mbox{and} \quad
\L_y = : [\, 0 \to e_n B^\tau \to V_y \to 0 \,].
\end{equation}
We can extend isomorphism $M_x \rTo^{e_n \kappa \cdot} M_y$ to
isomorphism of of the above complexes:
\[
\begin{diagram}[small, tight]
M_x                 & \rTo^{r_x}           & \L_x       \\
\dTo^{e_n \kappa \cdot} &                      & \dTo_{\Phi}\\
M_y                 & \rTo^{r_y}           & \L_y        \\
\end{diagram}
\]
First let us introduce some notations. Let $ B^\tau [S_1^{-1}] [S_2^{-1}]$
be a $\Gamma$-module where $B^\tau$ first localized by the set $S_1$
and next by $S_2$. Then it is easy to see that $e_n B^\tau [S_1^{-1}] [S_2^{-1}]
\cong e_n \mathbb C(x)(y)$  and $e_n B^\tau [S_2^{-1}] [S_1^{-1}]
\cong e_n \mathbb C(y)(x)$. We now introduce four linear maps:
\begin{equation}
\la{maps1}
\begin{diagram}[small, tight]
        &                 &e_n B^\tau [S_1^{-1}] [S_2^{-1}]   &
       &    \\
        & \ldTo^{\rrho_x} &   & \rdTo^{\lrho_y}   \\
e_n \mathbb C [x](y) &                 &   & & e_n \mathbb C (x)[y] \\
\end{diagram}
\qquad \quad
\begin{diagram}[small, tight]
        &                 & e_n B^\tau [S_2^{-1}] [S_1^{-1}]   &
        &    \\
        & \ldTo^{\lrho_x} &   & \rdTo^{\rrho_y}   \\
e_n \mathbb C (y)[x] &                 &   & & e_n \mathbb C [y](x) \\
\end{diagram}
\end{equation}
which are defined as follows $ \rrho_x : e_n f(x) g(y) \mapsto e_n f(x)_{+}
g(y)$, \,
$ \lrho_y : e_n f(x) g(y) \mapsto e_n f(x) g(y)_{+}$, \,
$\lrho_x : e_n g(y) f(x) \mapsto e_n g(y)f(x)_+$ and
$ \rrho_y : e_n g(y) f(x) \mapsto e_n g(y)_+f(x)$. It is clear that all of these maps are
$\Gamma$-equivariant. We then define a $\Gamma$-equivariant map
 $\phi : e_n B^\tau \to e_n B^\tau$ by
\begin{equation}
\la{phi}
\phi (e_n b) := \rrho_y\,\lrho_x(e_n \kappa \cdot e_n b) 
 = \rrho_y\,\lrho_x(e_n \kappa b)\ , \quad b \in
B^\tau\ .
\end{equation}
Now one can argue as in Lemma $7$ of \cite{BC} to prove the following
result.
\begin{proposition}
\la{lphi}
Let $\phi : e_n B^\tau \to e_n B^\tau$ be a map as in \eqref{phi} then:

$(1)\ $ $ \phi \,$ extends $\, \kappa \,$ through $\, r_x \,$,
i.~e. $\, \phi \circ r_x = r_y \circ e_n \kappa \,$.

$(2)\ $ $ \phi $ is invertible with
$\, \phi^{-1}: e_n B^\tau \to e_n B^\tau \,$ given
by $ \phi^{-1}(a) = \rrho_x\,\lrho_y (e_n \mu b) \,$.

$(3)\ $ We have $\, \phi(e_n b) = e_n b \,$ whenever
$\, b \in \mathbb C[x] $ or $ b \in \mathbb C[y]\,$.
\end{proposition}
\begin{proof}
Denote by $\, \mathbb C (x)_{-} \,$ the subspace of $ \mathbb C(x) $ consisting
 of functions vanishing at infinity. Then we can extend our earlier notation writing, 
for example, $\, \mathbb C (x)_{-}(y) \,$ for the subspace of $ \mathbb C (x)(y) $ 
spanned by all elements $\, f(x)\,g(y) \,$ with $\, f(x) \in k(x)_{-}\,$ and $\, g(y) \in k(y)\,$.

($1$)\,  Since $\, M_x \subset e_n \mathbb C(x)[y] \,$ we have $\, r_x(m) - m 
\in e_n \mathbb C(x)_{-}[y] = \mathbb C [y](x)_{-}\,$ for any $ m \in M_x\,$. 
Hence,$\, e_n \kappa \cdot ( r_x(m) - m) \in e_n \mathbb C(y)(x)_{-}\,$ 
and therefore $\, \lrho_x(e_n \kappa \cdot r_x(m)) = \lrho_x (e_n \kappa \cdot m)\,$.
On the other hand if $m\in M_x$ then $e_n \kappa \cdot m \in M_y \subset 
\mathbb C(y)[x]$ and therefore $\lrho_x(e_n \kappa \cdot m) = e_n \kappa \cdot m$.
Thus we have $\phi(r_x(m)) = \rrho_y \lrho_x (e_n \kappa \cdot m) =
\rrho_y (e_n \kappa \cdot m) = r_y (e_n \kappa \cdot m)$.

$(2)\,$  It follows from definition of $\phi$ that 
$\, \rrho_y\,\lrho_x(\phi(e_n b) - \kappa \cdot e_n b) = 0 \,$ and therefore 
$$\phi(e_n b) - e_n\kappa \cdot b \, \in \, e_n \mathbb C (y)(x)_{-} + e_n \mathbb C(y)_{-}(x)  = 
e_n \mathbb C(x)_{-}[y] + e_n \mathbb C[x](y)_{-} + e_n \mathbb C(y)_{-}(x)_{-}\ .
$$
Now multiplying the last expression by $ e_n \mu $ and using the fact that 
$\, e_n \mu - 1 \in e_n \mathbb C(x)_{-}(y)_{-} $, we obtain
$$
e_n \mu\cdot \phi(e_n b) - e_n b \, \in \, e_n \mathbb C(x)_{-}(y) + e_n \mathbb C(x)(y)_{-} + 
e_n \mathbb C(y)_{-}(x)_{-} + e_n \mathbb C(x)_{-}(y)_{-}(x)_{-}\ . 
$$
On the other hand, since $\, \phi(e_n b) \in e_n B^\tau \,$, we have
$\, e_n \mu \cdot \phi(e_n b) - e_n b \in e_n \mathbb C(x)(y) \,$.
Comparing the last two inclusions shows 
$\,e_n\mu \cdot \phi(e_n b) - e_n b \, \in \, e_n \mathbb C(x)_{-}(y) + e_n \mathbb C(x)(y)_{-} \,$. 
Hence $\, \rrho_x\,\lrho_y (e_n \mu \cdot \phi(e_n b) - e_n b) = 0 \,$ and 
therefore $\, \rrho_x\,\lrho_y (e_n \mu \cdot \phi(e_n b)) = e_n b \,$ for all 
$ b \in B^\tau \,$. Defining now $\, \phi^{-1}: e_n B^\tau \to e_n B^\tau \,$ by the formula 
$\, \phi^{-1}(e_n b) := \rrho_x\,\lrho_y (e_n \mu \cdot b) \,$ we see that
$ \phi^{-1} \circ \phi = \id_{e_n B^\tau} \,$. On the other hand, reversing the roles 
of $ \phi $ and $ \phi^{-1} $ in the above argument would give
obviously $\,\phi \circ \phi^{-1} = \id_{e_n B^\tau} \,$. Thus, $ \phi $ is an isomorphism 
of vector space, and $ \phi^{-1} $ is indeed its inverse.
                                                                                                
$(3)$ is immediate from the definition of $ \phi\,$. For example, if 
$ b \in \mathbb C[x] $  then $ e_n \kappa \cdot b - e_n b \in \mathbb C(y)_{-}(x) $ and therefore
$ \phi(e_n b) := \rrho_y\,\lrho_x(e_n \kappa \cdot b) = \rrho_y\,\lrho_x(e_n b) = b \,$,
as claimed.                               
\end{proof}
\begin{remark}
Once the isomorphism $\phi$ satisfying condition $(1)$ of
Proposition~\ref{lphi} is established one can easily determine
isomorphism of quotient spaces $\bar{\phi} : V_x \to V_y$ and hence
isomorphism of complexes $\Phi = (\phi,\, \bar{\phi}) : \L_x \to
\L_y$.
\end{remark}

We will now define our DG module. Let $\L := \L_x$ and endomorphisms
$X, \, Y \in End_{\mathbb C} (L^0)$ and $\X, \, \Y \in End_{\mathbb C}
(L^1)$
are given by
\begin{equation}
\la{xY1}
X(e_n b) : = \phi^{-1} (\phi (e_n b) \cdot x) , \quad Y(e_n b) =e_n b \cdot y
\end{equation}
\begin{equation}
\la{xY2}
\X(e_n b) : = \bar \phi^{-1} (\bar \phi (\overline{e_n b})
\cdot x) , \quad \Y(e_n b) =e_n b \cdot y
\end{equation}
where $"\cdot"$ stands for usual multiplication in $B^\tau$. It is
clear from the construction that these endomorphisms satisfy
\eqref{psi1}. We next define the 'cyclic' vectors:
\begin{equation}
\la{cvc}
i: W_n \to L^0 \, , \, e_n \mapsto e_n , \mbox{ and } \, \i :W_n \to L^1  
\, , \, e_n \mapsto d_{\L}(e_n).   
\end{equation}
Now the condition \eqref{psi2} is a consequence of the
following proposition.
\begin{proposition}
\la{P7}
The endomorphisms \eqref{xY1} and \eqref{xY2} satisfy the equations
\begin{equation}
\la{hghg}
X \, Y - Y\, X + T = i\, j \,  , \quad \X \, \Y - \Y\, \X + \bar T = \i \,\j
\end{equation}
for some $ j : L^0 \to W_n $ and $ \j : L^1 \to W_n $ related by  $ j = \j
d_{\L}$.
\end{proposition}
\begin{proof}
It suffices to show that 
$$ X \, Y(e_n b) - Y\, X(e_n b) + T(e_n b) \in \mathbb Ce_n \mbox{  for any } b\in B^\tau .$$
Indeed, if it holds we can define $j(e_n b) =X \, Y(e_n b) - Y\, X(e_n b) + T(e_n b)$.  
By previous Proposition it is then easy to see that 
$\, j(e_n b) = 0 \,$ on $\, \im(r_x) \,$, and since $ \im(r_x) = \Ker(d_{\L}) $
the second equation follows from the first.

Let $\tilde {b}:= X(e_n b) -e_n b \cdot x $, then using  \eqref{xY1} we have
$$  \phi (\tilde {b}) = \phi(e_n b)\cdot x - \phi(e_n b\cdot x)              
=  \rrho_y\,(\lrho_x(e_n \kappa b) \cdot x - \lrho_x(e_n \kappa b \cdot x)). $$
It is clear that the last expression lies in $e_n \mathbb C[y]$ and therefore,
by Proposition \ref{lphi}($3$) we get $\tilde {b} \in e_n \mathbb C[y]$ for all $b\in B^\tau$.
Now we have 
\begin{eqnarray}
\lefteqn{(XY - YX)(e_n b) + T(e_n b) = \phi^{-1}(\phi(e_n by)x) - \phi^{-1}(\phi(e_n b)x)y + T(e_n b) = } 
\nonumber \\
&& \qquad \left(\phi^{-1}(\phi(e_nb y)x) - e_n byx \right) - 
\left(\phi^{-1}(\phi(e_n b)x) - e_n bx \right)y\, \in \, e_n \mathbb C[y]\ .\nonumber
\end{eqnarray}
On the other hand, if besides (\ref{xY1}) we define $\, X',\,Y' \in \End_{\mathbb C}(e_n B^\tau)\,$ by
$$
X'(e_n b) := e_n b \cdot x \ ,\quad Y'(e_n b) := \phi(\phi^{-1}(e_n b) \cdot y)
$$
then by symmetry $\, (X'Y' - Y'X')e_n b + e_n b \in e_n \mathbb C[x] \,$ for all $ b \in B^\tau\,$. 
But $\,\phi \, X = X'\,\phi $ and $\,\phi \, Y = Y'\,\phi \,$.
Hence $\,\phi\left([X,\,Y]e_n b + e_n b\right) = [X',\,Y']\phi(e_n b) + \phi(e_n b) \in e_n \mathbb C[x]\,$,
and therefore
$$
[X,\,Y]e_n b + e_n b \, \in \, e_n \mathbb C[y]\,\cap \, \phi^{-1}(\mathbb C[x]) =  
\mathbb C[y]\,\cap \, \mathbb C[x] = \mathbb C\ .
$$
where the first equality holds due to Propostion \ref{lphi}$(3)$.
\end{proof}

Now if we choose $ f= i \, \j$ then
$d_{\boldsymbol C}(f) = (f_1,\, f_2) = (i\, j, \, \i\, \j)$ and
therefore, by Proposition~\ref{P7} condition   \eqref{psi2} holds.
\section{Bijective Correspondences}

Let us remind that $\mathcal R (V,W) $ is the set of isomorphism classes
of projective modules $M$  over $B^{\tau}$  such that
$[M]_{K_0}=[W] + [V]([L] - 2[W_0])$
under $K_0(B^{\tau}) \cong K_0(\Gamma)$.
Further let $\mathcal M(V,W)$ be the set of strict isomorphism classes of
DG-models as defined in Definition $1$.
Finally, let $\tilde {\mathfrak M}^{\tau}_{\Gamma}(V,W) =
\bigsqcup_{k=0}^{\infty}\mathfrak M^{\tau}_{\Gamma}
(V \oplus {\mathbb C}\Gamma ^{\oplus k},W)$
be a disjoint union of Nakajima spaces
  then we establish the following bijective correspondences.
\begin{theorem}
\la{T6}
There are four maps
\begin{equation}
\la{D7}
\begin{diagram}[small, tight]
\mathcal R(V,W) &  \pile{\rTo^{\theta_1}\\ \lTo_{\omega_1}} & \mathcal 
M(V,W) &
\pile{\rTo^{\theta_2}\\ \lTo_{\omega_2}} & \tilde {\mathfrak M}^{\tau}_{\Gamma}(V,W)
\ , \\
\end{diagram}
\end{equation}
such that $ (\theta_1, \omega_1) $ and $ (\theta_2, \omega_2) $
are pairs of mutually inverse
bijections
\end{theorem}
\begin{proof}
The map $\, \theta_1 $ is given by the construction
in Section~\ref{exist} which assigns
to an ideal $ M $ its  DG-model $ M \stackrel{r}{\to} \L $
( Theorem~\ref{existence}). Passing from
$M$ to isomorphic module produces DG-model quasi-isomorphic to $\L$
which by uniqueness theorem implies that they are DG-isomorphic and
therefore
this map is well-defined.

The map $\, \omega_1 $ is defined simply by taking cohomology
of DG-model which is by definition projective module of
$B^{\tau} $ such that $\phi_1^{-1}([M]_{K_0})=[W] + [V]([L] - 2[W_0])$.
Now it is clear that  $\, \omega_1 \circ \theta_1 = Id_{\mathcal R}$, while
$\,  \theta_1  \circ \omega_1 = Id_{\mathcal M}$ follows again
from uniqueness theorem.

In Section $2.2$ we have constructed Nakajima data from DG-model.
Since the action of $ \boldsymbol B $ commutes with DG-module isomorphism
we get well-defined map $\, \theta_2 $ from $ \mathcal M$
to  $ \tilde {\mathfrak M}^{\tau}_{\Gamma}(V,W) $.

In Section $2.3$ we have shown how to get DG-model from a point in
$\tilde {\mathfrak M}^{\tau}_{\Gamma}(V,W)$. Now if we replace
$ (\X, \Y, \i,\j) $ by  equivalent data $\, (g\X g^{-1}, \,g\Y g^{-1},\, g(\i),\,\j g^{-1})$,
where $g \in \GL(V\oplus\mathbb C \Gamma^{\oplus k})\,$,
then the functional $\lambda$ remains the same,
and hence so do the ideal $J$
and the $ R$-module $ L^0 $. On the other hand, the
differential $ d_{\L} $ gets changed to
$\, g\,d_{\L} \,$. As a result, we obtain an
DG-module $ \Tilde {\L} $ strictly isomorphic
to $ \L \,$, the isomorphism $ \L \to \Tilde {\L} $
being given by $\, (\id_{L^0}, g) \,$.
Thus, the  construction of Section $2.3$ yields a well-defined map
$ \omega_2: \tilde {\mathfrak M}^{\tau}_{\Gamma}(V,W) \to \MM\,$.

Now we have to show that $\,  \theta_2  \circ \omega_2 = Id$ and 
$\, \omega_2 \circ \theta_2 =Id_{\MM}$. The first equality follows immediately  from the
constructions in Sections $2.2$ and $2.3$. The second equality
follows from Theorem \ref{uniq} since both $\L$ and $\, \omega_2
\circ \theta_2 (\L)$ have the the same linear data $ (\X, \Y, \i,\j)
$ and hence produce the same $\lambda$.
\end{proof}
\section{ $G$ - equivariance}
Let $G=Aut_{\Gamma}(R)$ be the group of $\Gamma$-equivariant
automorphisms of the algebra $R=\mathbb C\langle x,y \rangle\ast
\Gamma$ preserving the form $\omega=xy-yx\in R$. In this section
we show that $G$ acts naturally on each of the spaces $\mathcal
R(V,W)$, $\mathfrak M(V,W)$ and $\mathcal M^{\tau}_{\Gamma}(V,W)$
and the bijections of Theorem \eqref{T6} are equivariant with
respect to these actions.

We start by describing the action of $G$ on the space of ideals
$\mathcal R(V,W)$. First, we observe that $G$ maps to the group
$Aut_{\Gamma}(B^{\tau})$ of $\Gamma$-equivariant automorphisms of
the algebra $B^\tau$ as $B^\tau$ is, by definition, a quotient of
the algebra $R$. Now, $Aut_{\Gamma}(B^\tau)$ acts naturally on the
category $\Mod(B^\tau)$ by twisting the structure of
$B^\tau$-modules by automorphisms: to be precise, for each
$\sigma\in Aut_{\Gamma}(B^\tau)$ we have an auto-equivalence
$\sigma_{\ast}: \Mod(B^\tau) \to \Mod(B^\tau)$, given by
$\sigma_{\ast}(M)=M_{\sigma^{-1}}$. Clearly, the functors
$\sigma_{\ast}$ restrict to the subcategory $\mathtt{PMod}
(B^\tau)$ of f.g. projective $B^\tau$-modules and their action
preserves the rank of projective modules. Thus, for each $\sigma
\in Aut_{\Gamma}(B^\tau)$ we have a bijection $\mathcal R \to
\mathcal R$ induced by $\sigma_{\ast}$, and this defines an action
of $G$ on $\mathcal R$ via the group homomorphism $G\to
Aut_{\Gamma}(B^\tau)$. We claim
\begin{lemma}
\la{gact}
 The action of $G$ on $\mathcal R$ defined above
respects the stratification $(2.11)$.
\end{lemma}
\begin{proof}
The action of the group $G$ on the category
$\mathtt{PMod}(B^\tau)$ by exact additive functors yields a
well-defined group homomorphism $ G \to Aut_{\Gamma}(B^\tau) \to
Aut(K_0(B^\tau))$; thus for each $\sigma \in G$, we have an
abelian group automorphism $\sigma_{\ast}: K_0(B^\tau) \to
K_0(B^\tau), [M]_{K_0} \mapsto [M_{\sigma^{-1}}]$. Now, in the
view of Lemma $2$, if $M\in \mathcal R(V,W)$, its stable
isomorphism class $[M]_{K_0}$ can be decomposed as
\begin{equation}
\la{asd}
 [M]_{K_0}=[W\otimes_{\mathbb C\Gamma}B^\tau]_{K_0}+[(V\otimes
L)\otimes_{\mathbb C\Gamma}B^\tau]_{K_0}-2[V\otimes_{\mathbb
C\Gamma}B^\tau]_{K_0}
\end{equation}
Since $\sigma \in G$ is $\Gamma$-equivariant, the corresponding
algebra automorphism $\sigma: B^\tau \to B^\tau$ an isomorphism
$B^\tau \cong (B^\tau)_{\sigma^{-1}}$ of $\mathbb
C\Gamma$-$B^\tau$-bimodules. Hence, with decomposition
\eqref{asd}, we set at once that
$[M_{\sigma^{-1}}]_{K_0}=[M]_{K_0}$ for every $M \in \mathcal R$
and $\sigma \in G$.This finishes the proof of the lemma.
\end{proof}

Thus, with Lemma \ref{gact}, we can define an action of the group
$G$ on $\mathcal R(V,W)$ simply by restricting its natural action
on $\mathcal R$.

Next, we define an action of $G$ on $\mathcal M(V,W)$. Again, we
start by observing that $G$ maps naturally to the group
$\mathtt{DGAut}_{\Gamma}(\boldsymbol B)$ of $\Gamma$-equivariant
automorphisms of the DG-algebra $\boldsymbol B$: in fact, given
$\sigma \in G$, we define $\tilde{\sigma}\in
\mathtt{DGAut}_{\Gamma}(\boldsymbol B)$ on generators by
$\tilde{\sigma}(x)=\sigma(x), \tilde{\sigma}(y)=\sigma(y),
\tilde{\sigma}(\nu)=\nu$. Each $\tilde{\sigma}\in
\mathtt{DGAut}_{\Gamma}(\boldsymbol B)$ yields an autoequivalence
$\tilde{\sigma_{\ast}}:\DGMod(\boldsymbol B) \to
\DGMod(\boldsymbol B)$ by twisting the action of $\boldsymbol B$
by $\tilde{\sigma}^{-1}$. It is clear that such autoequivalences
preserve the class of DG-models, since each axiom of Definition
$2$ is stable under twisting by $\tilde{\sigma}\in
\mathtt{DGAut}_{\Gamma}(\boldsymbol B)$. Moreover, if $\L\in
\mathcal M(V,W)$, then $H^0(\L) \in \mathcal R(V,W)$ and hence
$\sigma_{\ast}(H^0(\L))\in \mathcal R(V,W)\Longrightarrow
\tilde{\sigma}_{\ast}(\L) \in \mathcal M(V,W)$ by Lemma $6$. Thus,
the above action of $G$ on DG-models preserves each stratum
$\mathcal M(V,W)$, and it is obvious that the bijections
$\theta_1$ and $\omega_1$ are $G$-equivariant with respect to this
action and the action of $G$ on $\mathcal R(V,W)$ defined in Lemma
\ref{asd}.

Finally, it remains to define an action of $G$ on the quiver
varieties $\mathfrak M^{\tau}_{\Gamma}(V,W)$. To this end, as in
Section 2, we represent the points of $\mathfrak
M^{\tau}_{\Gamma}(V,W)$ by quadruples of matrices $(\X,\Y,\i,\j)$
and let
$\sigma.(\X,\Y,\i,\j):=(\sigma^{-1}(\X),\sigma^{-1}(\Y),\i,\j)$.
Since $\sigma$ is $\Gamma$-equivariant and preserves the form
$\omega=xy-yx$, this action is well-defined: the quadruple
$\sigma.(\X,\Y,\i,\j)$ satisfies the relations $(2.7)$ and
$(2.8)$. Moreover, it is clear that $\sigma.(\X,\Y,\i,\j)$ are
precisely the Nakajima data corresponding to the ``twisted''
DG-model $\tilde{\sigma_{\ast}}(\L)$ if $(\X,\Y,\i,\j)$
corresponds to $\L$. Thus, we have an action of $G$ on $\mathfrak
M^\tau_{\Gamma}(V,W)$ such that the bijection $\theta_2, \omega_2$
are $G$-equivariant. Summing up, we have estabilished the
following

\begin{theorem}
\la{T7}
The maps $ (\theta_1, \omega_1) $ and $ (\theta_2, \omega_2) $ are 
$G$-equivariant
bijective correspondences.
\end{theorem}

\section{Invariant Subrings of the Weyl Algebra}
In this section we look at the simpliest example of the algebra $O^\tau$ corresponding to $\tau=1$.
It is well-khown that in this case the algebra $B^\tau$ is isomorphic to the crossed product 
$A_1(\mathbb C)\ast \Gamma$ and $O^\tau$ to the subring $A_1^{\Gamma}$ of invariant of the
first Weyl algebra $A_1(\mathbb C)=\mathbb C \langle x,y \rangle / (xy-yx-1)$ under the action 
$x\mapsto \epsilon x$ and $y\mapsto \epsilon^{-1} y$. In fact, we have
$$ B^\tau \cong \mathbb C \langle x,y \rangle \ast \Gamma / (xy-yx-1) \cong
\left (\mathbb C \langle x,y \rangle / (xy-yx-1) \right ) \ast \Gamma \ ,$$ 
$$ O^\tau \cong A_1^{\Gamma}(\mathbb C) \ .$$

It is a simple observation that, in the case of $\tau=1$, the Nakajima variety can be embedded
as an affine algebraic subvariety in the Calogero-Moser variety. Indeed, from the relation \eqref{ma}
for a quadruple $(\X,\Y,\i,\j) \in \mathfrak M^\tau_{\Gamma}(V,W)$ we have
\begin{equation}
\la{ma23}
 \X\Y-\Y\X +I=\i\j \,
\end{equation}
which is exactly the Calogero-Moser relation. Now, a pair $(\X,\Y)$ satisfying this relation
does not have common invariant subspace (see Lemma $1.3$, \cite{W}) and hence,
the condition $(ii)$ of \eqref{NV} in the definition of the Nakajima variety 
$\mathfrak M^\tau_{\Gamma}(V,W)$ is redundant. Thus, in the case of $\tau=1$, the Nakajima 
variety is a subvariety of the Calogero-Moser variety whose points also satisfy the equations 
\begin{equation}
\la{mg23}
\X\,\bar{G} =\epsilon(g)\, \bar{G}\,\X  \quad  , \quad  \Y\,\bar{G}
=\epsilon^{-1}(g)\, \bar{G}\,\Y
\end{equation}
Now we will give another description of the Nakajima variety.
For this we remind that $\{W_0, W_1, ..., W_{m-1}\}$ is the complete set of irreducible 
$\Gamma$-modules such that the character of $W_i$ is $\epsilon^{i}$.
Then, if  $V \cong \bigoplus_{i=0}^{m-1} V_i \otimes W_i$ is the irreducible $\Gamma$-decomposition 
of $V$, we have
$$ Hom_{\Gamma}(V,V\otimes \epsilon) \cong \bigoplus_{i=0}^{m-1} Hom(V_i,V_{i-1})\, , \,
  Hom_{\Gamma}(V,V\otimes \epsilon^{-1}) \cong \bigoplus_{i=0}^{m-1} Hom(V_{i},V_{i+1}) $$
$$Hom_{\Gamma}(W_n,V) \cong  Hom(\mathbb C,V_n) \, , \,
Hom_{\Gamma}(V,W_n) \cong  Hom(V_n,\mathbb C)\, .$$

We now introduce the following algebraic variety (see \cite{N1}):
\begin{eqnarray}
&&\hspace{-32mm} D^n_{(k_0,...,k_{m-1})}:= \Bigl\{ \Bigl(\X_0, \X_1, ..., \X_{m-1}; \Y_0,\Y_1,..., \Y_{m-1}, \i_n, \j_n\Bigr)\,
\Big |  \\*[1ex] \nonumber
 && \hspace{-2mm}  \X_i \in Hom(V_{i+1}, V_i)\ , \ \Y_i \in Hom(V_{i},V_{i+1}) \, ,  \\*[1ex] \nonumber
 &&   \hspace{-2mm} \i_n \in Hom(\mathbb C, V_n) \, ,\, \j_n \in Hom(V_n, \mathbb C) , \\*[1ex] \nonumber
 &&  \hspace{-2mm} \X_i\Y_i - \Y_{i-1}\X_{i-1} + Id_{k_i} =0 \ , \, i\neq n , \\*[1ex] \nonumber
  &&   \hspace{-2mm} \X_n\Y_n - \Y_{n-1}\X_{n-1} + Id_{k_n} =\i_n\j_n  \Bigr\} \Big/\!\!\!\Big/ \prod_{i} GL(V_i) \, , \nonumber
\end{eqnarray}
where $k_i:=dim_{\mathbb C}(V_i)$. Then, due to equations \eqref{ma23} and \eqref{mg23}, there is a well-defined map 
$$\psi: \mathfrak M^{\tau}_{\Gamma}(V,W_n) \longrightarrow D^n_{(k_0,...,k_{m-1})} \ , $$ 
$$ \X \mapsto (\X_0,\X_1,...,\X_{m-1})\, , \, \Y \mapsto (\Y_0,\Y_1,...,\Y_{m-1}) \, ,
\, \i\mapsto \i_n \, , \, \j \mapsto \j_n \, .$$
 In fact one can easily prove the following result:
\begin{theorem}
The map $\psi$ is an isomorphism of algebraic varieties with inverse map defined by 
$$ (\X_0,\X_1,...,\X_{m-1}) \mapsto \X \, , \, (\Y_0,\Y_1,...,\Y_{m-1}) \mapsto \Y  \, ,
\, \i_n \mapsto \i \, , \,  \j_n \mapsto \j \ , $$
where $\X$ and $\Y$ are the following matrices
\begin{equation}
\X=
\begin{pmatrix} 
0 & \X_0 & 0& \hdots  &0 \\
 0& 0 & \X_1  & \hdots & 0\\
0& 0&0&\ddots&\vdots \\
\vdots& \vdots&\ddots&\ddots &\X_{m-2}\\
\X_{m-1} &0  & \hdots &0& 0 
\end{pmatrix} 
\, , \,
\Y=
\begin{pmatrix} 
0 & 0 &0& \hdots &  \Y_{m-1} \\
\Y_0 & 0 &0& \hdots &  0\\
0 & \Y_1 &0& \ddots &\vdots \\
\vdots& \vdots&\ddots& \ddots &0\\
0 &0 & \hdots & \Y_{m-2} & 0
\end{pmatrix} 
 \, .
\end{equation}
Moreover, we have 
\begin{equation}
\la{dimnv1}
 dim_{\mathbb C}\, D^n_{(k_0, k_1)}=
2(k_n -(k_0-k_1)^{2}) \ , \mbox{ for } m=2 \ ,
\end{equation}
\begin{equation}
\la{dimnv}
 dim_{\mathbb C}\, D^n_{(k_0,...,k_{m-1})}=
2\bigg(k_n -\bigg(\sum_{i=0}^{m-1} k^2_i -\sum_{ i < j} k_i k_j\bigg)\bigg) \ , \mbox{ for } m>2 .
\end{equation}
\end{theorem}

Let us define the following set 
$$N_n=\bigg\{(k_0,k_1,...,k_{m-1}) \in \mathbb N^{m} \, | \, 
k_n -\bigg(\sum_{i=0}^{m-1} k^2_i -\sum_{ i < j} k_i k_j \bigg ) \ge 0\,\bigg\} \ .$$
From \eqref{dimnv} we can see that this set consists exactly of those points of 
$\mathbb N^{m}$ for which the corresponding Nakajima variety
$D^n_{(k_0,...,k_{m-1})}$ is nonempty. 
With this notation we can restate Theorem \ref{bgk} as follows
\begin{corollary}
\la{c5}
The set $\mathcal R(A_1^{\Gamma})$ of isomorphism classes of ideals of $A_1^{\Gamma}$ 
is in bijection with the union of algebraic varieties
$$\bigsqcup_{n=0}^{m-1}
\bigsqcup_{(k_0,...,k_{m-1})\in N_n} D^n_{(k_0,...,k_{m-1})} \, .$$
\end{corollary}
In the case $\,m=2\,$, the varieties $\,D^n_{(k_0,...,k_{m-1})}\,$ have been introduced 
recently in \cite{DNM} (see {\it loc. cit.}, Theorem~$3$) to classify the ideals of the  ${\mathbb Z}_2$-invariant 
subring of $A_1(\C) $. Our Corollary~\ref{c5} may be viewed thus as a generalization of this description to the 
case of an arbitrary cyclic group $\, \Z_m $.
\section{Appendix: $\A$-morphisms of DG modules}
The DG-algebra $\boldsymbol B$ regarded as  $\A$-algebra has only two structure
maps
$ m^{\boldsymbol B}_{1} := d_{\boldsymbol B}$ and  $m^{\boldsymbol B}_{2}$ the usual
associative multiplication in $\boldsymbol B$. For any DG-module $\LB$
over $\boldsymbol B$ viewed as $\A$-module has $m^{\LB}_{n} = 0$ for
$n \geq 3$ which satisfy Leibnitz rule:
$$ m^{\LB}_{1}  m^{\LB}_{2} =  m^{\LB}_{2} (m^{\LB}_{1}\otimes 1)
+  m^{\LB}_{2} (1\otimes d_{\boldsymbol B}) $$
Now we remind the definition of morphisms of $\A$-modules (see \cite{K})
\begin{definition}
\la{D2}
A {\it morphism of $\A$-modules} $ f : \LB \to \E $ is a sequence of graded
morphisms
\begin{equation}
\la{morph1}
f_{n} : \LB\otimes {\boldsymbol B}^{\otimes n-1} \to \E
\end{equation}
of degree $1-n$ such that for each $n \geq 1$, we have
\begin{equation}
\la{morph2}
\sum (-1)^{r+st} f_{u} \circ (1^{\otimes r} \otimes
m_{s} \otimes 1^{\otimes t}) =
\sum (-1)^{(r+1)s} m_{u} \circ (f_{r} \otimes 1^{\otimes s}) ,
\end{equation}
where the left hand sum is taken over all decompositions $ n =
r+s+t, r, t \geq 0, s\geq 1$ and we put $u=r+1+t$; and the right
hand sum is taken over all decompositions $n = r+s, r\geq 1, s\geq 0
$ and we put  $ u = 1+s $
\end{definition}
\begin{lemma}
\la{L3}
Let $ \LB $ and $\E$ be DG-modules over $\boldsymbol B$, $\LB$ having nonzero
components only in degree $0$ and $1$ and $\E$ positively graded:
$ \LB = L^0 \oplus L^1 $ and $\E = E^0 \oplus E^1 \oplus E^2 ...$.

$(a)$ Any $\A$-morphism $ f : \LB \to \E $ is determined by two
components $ (f_{1},f_{2}) $ satisfying the relations:
\begin{equation}
\la{3.3}
m^{\E}_{1} f^0_{1} = f^1_{1} m^{\LB}_{1}
\end{equation}
\begin{equation}
\la{3.4}
f^0_{1} (m^{\LB}_{2} (u, a)) - m^{\E}_{2}(f^0_{1}(u), a) =
f_{2}(m^{\LB}_{1}(u), a),
\quad  \forall\, u \in L^0, a \in R
\end {equation}
\begin{equation}
\la{3.5}
f^1_{1} (m^{\LB}_{2} (v, a)) - m^{\E}_{2}(f^1_{1}(v), a) =
m^{\E}_{1}f_{2}(v, a),
\quad  \forall\, v \in L^1, a \in R
\end {equation}
\begin{equation}
\la{3.6}
f^0_{1} (m^{\LB}_{2} (v, c)) - m^{\E}_{2}(f^1_{1}(v), c) = -f_{2}(v,
d_{\boldsymbol B}c),
\quad  \forall\, v \in L^1, c \in I
\end {equation}
\begin{equation}
\la{3.7}
f_{2} (v, ab) = m^{\E}_{2}(f_{2}(v, a), b) + f_{2}(m^{\LB}_{2}(v, a), b),
\quad  \forall\, u \in L^0,  a, b \in R
\end {equation}

$(b)$ if $m^{\LB}_{1}$ is surjective then equations \eqref{3.5}-\eqref{3.7}
are
formal consequences of \eqref{3.3} and \eqref{3.4}
\end{lemma}
\begin{proof}
The relation \eqref{3.3}  follows easily from \eqref{morph2} for $n =1$. For
$n =2$ we get the equation
\begin{eqnarray}
\la{3.8}
\lefteqn{-f_{2}(1 \otimes d_{\boldsymbol B}) + f_{1} \circ m^{\LB}_{2}
- f_{2}(m^{\LB}_{1} \otimes 1) =} \\*[1ex]
&&m^{\E}_{2}(f_{1} \otimes 1) + m^{\E}_{1} \circ f_{2}.   \nonumber
\end{eqnarray}
Since $ deg(f_{2}) = -1$ it has only one component
$f_{2}: L^1 \otimes R \to E^0 $ and
therefore the relations \eqref{3.4}-\eqref{3.6} are consequences
of \eqref{3.8}. For $ n=3 $ equation
\eqref{morph2} has the following form
\begin{eqnarray}
\la{3.9}
\lefteqn{f_{3}(1 \otimes 1 \otimes d_{\boldsymbol B}) +
f_{3}(1 \otimes d_{\boldsymbol B} \otimes 1)
+ f_{3}(m^{\LB}_{1} \otimes 1\otimes 1) -} \\
&& f_{2}(1\otimes m^{\boldsymbol B}_{2} ) + f_{2}(m^{\LB}_{2}
\otimes 1) + f_{1} \circ m^{\LB}_{3} =  \nonumber \\
&& m^{\E}_{3}(f_{1} \otimes 1\otimes 1) - m^{\E}_{2}
(f_{2} \otimes 1) + m^{\E}_{1} \circ f_{3}. \nonumber
\end{eqnarray}
By degree argument we can conclude that $ f_{n} =0$
for $ n \geq 3$. Now since both $\LB$ and $\E$
are DG-modules we have $ m^{\LB}_{3} = m^{\E}_{3} =0$.
The equation \eqref{3.9} can be simplified
\begin{equation}
\la{3.10}
- f_{2}(1\otimes m^{\boldsymbol B}_{2} ) + f_{2}(m^{\LB}_{2} \otimes 1) =
m^{\E}_{2}(f_{2} \otimes 1)
\end{equation}
which is equivalent to \eqref{3.6}-\eqref{3.7}.

To prove part $b$ we first apply $m_{1}$ to the equation
\eqref{3.4}. Then using \eqref{3.3} and
and $R$-linearity of $m_{1}$
(i.e. $d_{\boldsymbol B}(a) =0$ for any $a\in R $) we have
\begin{eqnarray}
&&m_{1}(f^0_{1} (m^{\LB}_{2} (u, a))) - m_{1}(m^{\E}_{2}(f^0_{1}(u), a)) =
\nonumber \\
&& f^1_{1} (m^{\LB}_{2} (m_{1}(u), a)) - m^{\E}_{2}(f^1_{1}(m_{1}(u)), a) =
m^{\E}_{1}f_{2}(m_{1}(u), a). \nonumber
\end{eqnarray}
Since $m_{1}$ is surjective this implies \eqref{3.5}.
Now let $a=d_{\LB}c$ and $ v =m^{\LB}_{1}(u)$
in \eqref{3.4} then
\begin{equation}
\la{3.61}
f^0_{1} (m^{\LB}_{2} (u, d_{\boldsymbol B}c)) - m^{\E}_{2}
(f^0_{1}(u), d_{\boldsymbol B}c) = f_{2}(m^{\LB}_{1}(u), d_{\boldsymbol B}c)
\end{equation}
Since $m^{\LB}_{2}(u, c) =0 $ for all $c \in I$
and $u \in L^0 $ by Leibnitz rule we get
$ m^{\LB}_{2}(m^{\LB}_{1}(u), c) = -m^{\LB}_{2}(u, d_{\boldsymbol B}c)$ .
Similarly one can show that
$m^{\E}_{2}(f^0_{1}(u), d_{\boldsymbol B}c) =-m^{\E}_{2}(f^1_{1}(v), c)$.
By plugging  the last two relations into  \eqref{3.61} we obtain
\eqref{3.6}.

Now we will show that \eqref{3.4} implies \eqref{3.7}. Let $v =m^{\LB}(u)$
then from
\eqref{3.4} we have
\begin{equation}
\la{3.12}
f_{2}(v, ab) = f_{2}(m^{\LB}_{1}(u), ab) = f^0_{1}(m^{\LB}_{2}(u, ab)) -
m^{\E}_{2}(f^0_{1}(u), ab)
\end{equation}
\begin{equation}
\la{3.13}
m^{\E}_{2}(f_{2}(v, a),b) = m^{\E}_{2}(f_{2}(m^{\LB}_{1}(u), a), b)
  = m^{\E}_{2}(f^0_{1}(m^{\LB}_{2}(u, a)), b) -
m^{\E}_{2}(f^0_{1}(u), ab)
\end{equation}
\begin{equation}
\la{3.14}
f_{2}(m^{\LB}(v, a), b) = f_{2}(m^{\LB}_{2}(m^{\LB}_{1}(u), a), b)
  = f^0_{1}(m^{\LB}_{2}(u, ab)) -
m^{\E}_{2}(f^0_{1}(m^{\LB}_{2}(u, a)), b).
\end{equation}
Adding now \eqref{3.13} and \eqref{3.14} and using
\eqref{3.12} we easily derive \eqref{3.7}.
\end{proof}
\end{document}